\documentclass{amsart}
\usepackage{graphicx}
\usepackage{amssymb}
\usepackage{epstopdf}
\usepackage{nicefrac}
\usepackage{enumerate}
\usepackage{verbatim}
\usepackage{hyperref}
\usepackage{color}
\usepackage{pst-all}
\usepackage{tikz}

\newtheorem{thm}{THEOREM}[section]

\newtheorem{defn}[thm]{DEFINITION}
\newtheorem{ex}[thm]{EXAMPLE}
\newtheorem{lemma}[thm]{LEMMA}

\newtheorem{prop}[thm]{PROPOSITION}
\newtheorem{quest}[thm]{QUESTION}
\newtheorem{remark}[thm]{REMARK}

\newcommand{\ds}{\displaystyle}
\newcommand{\mC}{{\mathbb C}}
\newcommand{\mQ}{{\mathbb Q}}
\newcommand{\mZ}{{\mathbb Z}}

\newcommand{\cA}{{\mathcal A}}
\newcommand{\cB}{{\mathcal B}}
\newcommand{\cC}{{\mathcal C}}
\newcommand{\cD}{{\mathcal D}}

\newcommand{\cM}{{\mathcal M}}

\newcommand{\cO}{{\mathcal O}}
\newcommand{\cP}{{\mathcal P}}

\newcommand{\oK}{{\overline{K}}}

 \newcommand\tail{\mathrel{\stackrel{\makebox[0pt]{\mbox{\normalfont\tiny t}}}{\sim}}}

\input xy
\xyoption {all}

\begin{document}

\title{Arboreal Cantor actions}

\begin{abstract}
{ In this paper, we consider minimal equicontinuous actions of discrete countably generated groups on Cantor sets, obtained from the arboreal representations of absolute Galois groups of fields. In particular, we study the \emph{asymptotic discriminant} of these actions.

 The asymptotic discriminant is an invariant 
 obtained by restricting the action to a sequence of nested clopen sets, and studying the isotropies of the enveloping group actions in such restricted systems. An enveloping (Ellis) group of such an action is a profinite group.

A large class of actions of profinite groups on Cantor sets is given by arboreal representations of absolute Galois groups of fields. 
We show how to associate to an arboreal representation an action of a discrete group, and give examples of arboreal representations with stable and wild asymptotic discriminant.

}
\end{abstract}

 \author{Olga Lukina}
 \email{lukina@uic.edu}
\address{Department of Mathematics, University of Illinois at Chicago, 322 SEO (m/c 249), 851 S. Morgan Street, Chicago, IL 60607-7045}

 \date{}


\thanks{Version date: August 30, 2018}

\date{}


\maketitle

\section{Introduction}\label{sec-intro}

Let $G$ be a countably generated discrete group, and let $X$ be a Cantor set, that is, a totally disconnected compact metrizable space, which has no isolated points (such a space is called \emph{perfect}). In this paper, Cantor sets arise as spaces of paths in $d$-ary rooted trees $T$, see Example \ref{ex-treecylinder} for more details. Such a space is also called the \emph{boundary} of the tree $T$ in geometric group theory \cite{Nekr}, or the set of \emph{ends} of a tree \cite{JonesLevy2017,Nekr}.

Let the group $G$ act on $X$ via the homomorphism $\Phi: G \to {Homeo}(X)$. We also denote the action by $(X,G,\Phi)$, and write $g\cdot x$ for $g(x)$. A group action $(X,G,\Phi)$ is \emph{minimal} if for every point $x \in X$, the orbit $\{g \cdot x \mid g \in G\}$ is dense in $X$. Let $D$ be a metric on $X$, compatible with the topology of $X$. An action $(X,G,\Phi)$ is \emph{equicontinuous} if for every $\epsilon >0$ there exists $\delta >0$ such that for every $g \in G$ and $x,y \in X$ with $D(x,y) < \delta$ we have $D(g\cdot x, g \cdot y) < \epsilon$. Intuitively, this definition means that for any two points $x,y \in X$, their images under the action of any element $g \in G$ are always at approximately the same distance from each other, never getting too far apart or too close together. In what follows, a `minimal Cantor system' or a `minimal Cantor action' or sometimes just `minimal action' always refers to a minimal action of a countably generated discrete group on a Cantor set.

When studying group actions on Cantor sets from the dynamical point of view, one would like to find invariants, which classify the actions up to a certain equivalence. For example, the most significant progress has been achieved in recent years in the classification of minimal actions of abelian groups up to orbit equivalence \cite{GPS1995,GMPS2008,GMPS2010}, and recently a few authors considered the classification of minimal equicontinuous actions of abelian groups \cite{GPS2017}, or, slightly more generally, free actions of non-abelian groups \cite{CM2016}, up to a conjugacy. 

The main interest of the current paper is equicontinuous minimal group actions on Cantor sets which have non-trivial isotropy groups, that is, there are elements in $G$ which fix some, but not all, points in the Cantor set $X$. Since the actions are minimal, points fixed by some elements in $G$ are necessarily moved around by other elements, as their orbits must be dense in $X$. In the setting of actions with non-trivial isotropy, the methods used for abelian or free actions need not work; for example, objects like the orbit cocycle in \cite{CM2016} may not be defined. Thus different, new invariants are needed to understand actions with non-trivial isotropy.

One such invariant, called the \emph{asymptotic discriminant}, has been introduced by the author joint with Hurder in \cite{HL2017a}, as a culmination of a series of papers on actions with non-trivial isotropy groups joint with Dyer and Hurder \cite{DHL2016a, DHL2016b,DHL2016c}. The asymptotic discriminant of a minimal equicontinuous action $(X,G,\Phi)$ is an invariant of return equivalence of minimal Cantor systems. Return equivalence for a more general class of pseudogroup actions was defined and is explained in more detail in \cite{CHL2017}. Recent work \cite{HL2017b} extends some methods from \cite{CM2016} to equicontinuous minimal Cantor systems with stable asymptotic discriminant.

In this paper, we show that a large class of equicontinuous minimal actions with non-trivial isotropy comes from so-called `arboreal representations' of absolute Galois groups into the groups of automorphisms $Aut(T)$ of $d$-ary rooted trees. The image of such a representation is a profinite group, acting on the tree $T$ by permuting infinite paths in the tree. As explained further in this section, the representation is obtained as the inverse limit of finite groups, and so it contains a countably generated dense subgroup $\Gamma$ \cite{RZ}. Denoting by $G$ the group $\Gamma$ with discrete topology, we obtain an action of a countably generated discrete group $G$ on the Cantor set of infinite paths in the tree.
This is made precise in our first main result, Theorem \ref{thm-suspension}, which is stated further in the introduction and proved in Section \ref{sections-provetwotheorems}. 

\begin{ex}\label{ex-treecylinder}
{\rm For readers with non-topological background we explain why the set of infinite paths in a tree $T$ is a Cantor set, and why the action of $G$ on the set of paths is minimal and equicontinuous. Recall that a Cantor set is a compact totally disconnected perfect metrizable space.

Let $d \geq 2$, and let $T$ be a $d$-ary rooted tree, that is, for $n \geq 0$ the set $V_n$ contains $d^n$ vertices, and every vertex in $V_n$ is connected by edges to precisely $d$ vertices in $V_{n+1}$. A path in $T$ is an infinite sequence $(v_n)_{n \geq 0} = (v_0,v_1,v_2,\cdots)$ such that $v_{n}$ and $v_{n+1}$ are connected by an edge, for $n \geq 0$. The set of all such sequences, which we denote by $\cP_d$, is a subset of the product $\prod_{n \geq 0} V_n$. The sets $V_n$ are given discrete topology, and so they are compact since they are finite. The product space $\prod_{n \geq 0} V_n$ is compact by the Tychonoff theorem \cite{Willard}, and $\cP_d$ is closed in $\prod_{n \geq 0}V_n$ by a standard argument. Thus $\cP_d$ is compact. Points are the only connected components in $\prod_{n \geq 0} V_n$, and so $\prod_{n \geq 0} V_n$ (and $\cP_d$) is totally disconnected.

Let $\widetilde{\pi}_n: \prod_{n \geq 0} V_n \to V_n$ be the projection, and denote by $\pi_n = \widetilde{\pi}_n|_{\cP_d}$ its restriction to $\cP_d$. Open sets in the product topology on $\prod_{n \geq 0} V_n$ have the form $\prod_{n \geq 0} U_n$, where $U_n \subseteq V_n$, and $U_n = V_n$ for all but a finite number of $n$. Let $U$ be such an open set, and let $m = \max \{n \mid U_n \ne V_n\}$. Consider the intersection $U \cap \cP_d$, which is an open set in $\cP_d$. Let $w_m \in U_m$ be a vertex, and note that $w_m \in U_{m}$ is joined by an edge to precisely one vertex $w_{m-1} \in V_{m-1}$. Inductively we obtain the vertices $w_i \in V_i$, $0 \leq i < m$, such that $w_{i+1}$ is joined by an edge to $w_i$. Since $\cP_d$ contains connected paths in the tree $T$, a sequence $(v_n)_{n \geq 0} \in \cP_d$ satisfies $v_m = w_m$ if and only if it satisfies the condition $v_i = w_i$ for $0 \leq i \leq m$. Thus the set of paths $U_m(w_m) = \{(v_n)_{n \geq 0} \mid v_m = w_m\} \subset \cP_d$ containing the vertex $w_m$, can be written as $U_m(w_m) = \pi_m^{-1}(w_m)$, and the open set $U \cap \cP_d$ can be written as $U \cap \cP_d = \pi_m^{-1}(U_m)= \cup\{\pi_m^{-1}(w) \mid w \in U_m\}$. See also the discussion before Proposition 2.2 in \cite{JonesLevy2017} for a description of $\cP_d$ using inverse limits.

Given a path $(v_n)_{n \geq 0}$, every open neighborhood of $(v_n)_{n\geq 0}$ contains an open set of the form $U_n(v_n)$ for some $n \geq 0$. Since $d \geq 2$, every such set $U_n(v_n)$ is infinite, and so $(v_n)_{n\geq 0}$ is not isolated. So $\cP_d$ is perfect. Since $V_n$ are finite sets, the complement of every open set $\pi_n^{-1}(U_n)$, where $U_n \subset V_n$, is also open, which means that $\pi_n^{-1}(U_n)$ is open and closed. A set which is open and closed is called a \emph{clopen} set.

Let $G$ be a countably generated discrete group, and let $G$ act on the tree $T$ by permuting vertices in each $V_n$, $n \geq 0$, in such a way that the connectedness of paths in $T$ is preserved, and the action is transitive on each $V_n$. Since permutations are bijective, the action of each $g \in G$ induces a bijective map $\Phi(g):\cP_d \to \cP_d$. The image of an open set $U_n(w)$ under $\Phi(g)$ is an open set $U_n(g \cdot w)$, so $\Phi(g)$ is a homeomorphism. Thus $G$ acts on $\cP_d$ by homeomorphisms.

Let $(v_n)_{n \geq 0} \in \cP_d$ be a path. Since $G$ acts transitively on $V_m$, for $m \geq 0$, for every $w_m \in V_m$ there exists $g \in G$ such that $g \cdot v_m = w_m$. Then $g \cdot (v_n)_{n\geq 0} \in U_m(w_m)$, and the orbit of $(v_n)_{n \geq 0}$ is dense in $\cP_d$. Thus $G$ acts minimally on $\cP_d$. 

To conclude Example \ref{ex-treecylinder}, define a metric $D$ on $\cP_d$ by setting
  \begin{align}\label{eq-metric} D((v_n)_{n \geq 0},(u_n)_{n \geq 0}) = \frac{1}{d^m}, \textrm{ where }m = \max \{ n \mid v_n = u_n\}, \end{align}
that is, $D$ measures the length of the longest finite path contained in both $(v_n)_{n \geq 0}$ and $(u_n)_{n \geq 0}$. 
Since $G$ acts bijectively on $V_n$ for $n \geq 0$, $(v_n)_{n \geq 0}$ and $(u_n)_{n \geq 0}$ contain a common path of length $m$ if and only if the images $g \cdot (v_n)_{n \geq 0}$ and $g \cdot (u_n)_{n \geq 0}$ contain a common path of length $m$, so the action of $G$ on $\cP_d$ is equicontinuous with respect to the metric $D$,  where we can take $\delta = \epsilon$ for every $\epsilon >0$.
}
\end{ex}

Many actions of discrete groups arising from arboreal representations of Galois groups have non-trivial isotropy; in fact, the actions with `as much isotropy as possible' are often of interest since they are related to certain questions in number theory \cite{Jones2014}. For example, an open problem in the area of arboreal representations is to determine for which rational functions the image of the arboreal representation has finite index in $Aut(T)$. In Theorem \ref{thm-wild} of this paper, for representations which have finite index in $Aut(T)$ we show that the associated discrete group action has non-stable (\emph{wild}) asymptotic discriminant. We give an example of an arboreal representation where the asymptotic discriminant is stable in Theorem \ref{thm-kummerextame}. In this theorem, the arboreal representation is given by the polynomial $f(x) = (x + p)^d - p$ over finite extensions of the $p$-adic numbers.

\medskip

We now prepare to state our main theorems, which will be proved in the subsequent sections. In order to do that we first need to review some background concepts and previous results.

Consider the action $\Phi: G \to { Homeo}(X)$ of a discrete group $G$ on a Cantor set $X$, and let $\widetilde{g}$ be the element in $Maps(X,X) = X^X$ defined by $g \in G$. The closure $E(X,G,\Phi)$ of the set $\{\widetilde{g} \, \mid \, g \in G\}$ in the topology of pointwise convergence is a compact subset of $X^X$ and has the structure of a semigroup, called the \emph{enveloping (or Ellis) semigroup}  \cite{Auslander1988,Ellis1969,Ellis2014}. If the action $(X,G,\Phi)$ is equicontinuous, then the Ellis semigroup $E(X,G,\Phi)$ is a group of homeomorphisms, and $E(X,G,\Phi)$ coincides with the closure $\overline{\Phi(G)}$ in $Homeo(X)$ in the uniform topology. For equicontinuous actions on Cantor sets the Ellis group $\overline{\Phi(G)}$ is a compact Hausdorff totally disconnected group. By compactness, open subgroups of $\overline{\Phi(G)}$ have finite index in $\overline{\Phi(G)}$, and so by \cite[Theorem 2.1.3]{RZ} $\overline{\Phi(G)}$ is a profinite group. The group $\overline{\Phi(G)} \subset Homeo(X)$ acts on the Cantor set $X$ by homeomorphisms. 

If the action $(X,G,\Phi)$ is minimal, that is, the orbits of points in $X$ are dense, then the action of $\overline{\Phi(G)}$ on $X$ is transitive, that is, the Cantor set $X$ is a single orbit of this action. For a point $x \in X$, denote by $\overline{\Phi(G)}_x = \{\widetilde{g} \in \overline{\Phi(G)} \, \mid \, \widetilde{g}\cdot x = x\}$ the isotropy subgroup of the Ellis group action at $x$. The subgroup $\overline{\Phi(G)}_x$ is closed, and the natural action of $G$ on the coset space $\overline{\Phi(G)}/\overline{\Phi(G)}_x$ is conjugate with $(X,G,\Phi)$ \cite{Auslander1988}.

Theorem \ref{thm-groupchains}, stated below and proved in\ \cite{DHL2016a}, allows us to explicitly compute the Ellis group and the isotropy subgroup of its action at a given point $x \in X$. For a given descending chain 
  $$\{G_n\}_{n \geq 0} = \{ G_0 \supset G_1 \supset \cdots \}, \, G_0 = G,$$ 
 of finite index subgroups of $G$, denote by ${\ds G_\infty = \lim_{\longleftarrow}\{G/G_{n+1} \to G/G_{n}\}}$ the inverse limit of coset spaces, with bonding maps induced by the inclusions $g G_{n+1} \subset g G_{n}$. There is a natural left action of $G$ on $G_\infty$, given by the left multiplication
  \begin{align}\label{eq-leftaction}g \cdot (g_0 G_0, g_1 G_1,\ldots) = (gg_0G_0, gg_1G_0, \ldots). \end{align}
We denote such action by $(G_\infty,G)$. An adaptation of the result in \cite{ClarkHurder2013} to group actions, which can be found in the Appendix in \cite{DHL2016a}, shows that if $(X,G,\Phi)$ is minimal and equicontinuous, and $x \in X$, then there exists a descending chain of finite index subgroups $\{G_n\}_{n \geq 0}$ and a homeomorphism  $\phi: X \to G_\infty$, such that $\phi(x) = {\bf e}$, where ${\bf e} = (eG_n)$ denotes the sequence of cosets of the identity $ e \in G$, and $\phi(g \cdot y) = g \cdot \phi(y)$ for all $y \in X$ and all $g \in G$. That is, $(X,G,\Phi)$ and $(G_\infty, G)$ are pointed conjugate.

For each $G_n$ in the group chain, denote by ${\ds C_n = \bigcap_{g \in G} gG_ng^{-1} }$ its maximal normal subgroup. Maximal normal subgroups form a descending chain $\{C_n\}_{n \geq 0}$, and their inverse limit 
  $${\ds C_\infty = \lim_{\longleftarrow}\{G/C_{n+1} \to G/C_{n}\}}$$ 
is a profinite group. The subgroup
  $${\ds \cD_x =  \lim_{\longleftarrow}\{G_{n+1}/C_{n+1} \to G_{n}/C_{n}\}}$$ 
contains elements of $C_\infty$ which fix $x$.

\begin{thm}\label{thm-groupchains}\cite{DHL2016a}
Let $(X,G,\Phi)$ be a minimal equicontinuous group action, and let $x \in X$. Let  $\{G_n\}_{n \geq 0}$ be a descending chain of finite index subgroups such that $(G_\infty, G)$ and $(X,G,\Phi)$ are conjugate. Then there exists an isomorphism $\widetilde{\phi}: \overline{\Phi(G)} \to C_\infty$ of the Ellis group onto the inverse limit group $C_\infty$, which restricts to the isomorphism $\widetilde{\phi}: \overline{\Phi(G)}_x \to  \cD_x$ of the isotropy groups.
\end{thm}

Although a standing assumption in \cite{DHL2016a,DHL2016c,HL2017a} was that $G$ is a finitely generated group, the reason for that was not any restrictions imposed by proofs or by the properties of the objects considered. The motivation in \cite{DHL2016a,DHL2016c,HL2017a} was to study and classify the dynamics of weak solenoids, and for a group to act on the Cantor fibre of a weak solenoid it must be realizable as a homomorphic image of a fundamental group of a closed manifold. This leads to the assumption of finite generation in  \cite{DHL2016a,DHL2016c,HL2017a}. However, the notion of the Ellis group does not require finite generation, and finite generation was not used in any of the proofs in \cite{DHL2016a,DHL2016c,HL2017a}. One easily checks that Theorem \ref{thm-groupchains} and all results in \cite{DHL2016a,DHL2016c,HL2017a} which we use in this paper are true for countably generated groups as well. However, some related results on strong quasi-analyticity in \cite{ALC2009,ALM2016} may require the finite (more precisely, compact) generation assumption.

The homeomorphisms $\phi$ and $\widetilde{\phi}$ in Theorem \ref{thm-groupchains} depend on the choice of a point $x \in X$ and a group chain $\{G_n\}_{n \geq 0}$, as discussed in more detail in \cite{DHL2016a}. The Ellis group $\overline{\Phi(G)}$ depends only on the action $(X,G,\Phi)$, and the isotropy group $\overline{\Phi(G)}_x$ depends on the choice of $x \in X$ up to an isomorphism. Thus the cardinality of the isotropy group $\cD_x$ is an invariant of the action $(X,G,\Phi)$.

The relationship between the cardinality of $\overline{\Phi(G)}_x$ and the properties of the action was studied in \cite{DHL2016a,DHL2016b,DHL2016c}. In particular, the isotropy group $\overline{\Phi(G)}_x$ is trivial if and only if the subgroups in the group chain $\{G_n\}_{n \geq 0}$ can be chosen to be normal (such actions are called $G$-odometers in \cite{CP2008}). For example, for the actions of abelian groups the isotropy group $\overline{\Phi(G)}_x$ is always trivial. Automorphisms of the minimal equicontinuous Cantor system $(G_\infty, G)$ (where $G$ acts on the left) are given by the right action of elements of $C_\infty$ on $G_\infty$. It is shown in \cite{DHL2016a} that the group of automorphisms acts transitively on $G_\infty$ if and only if the isotropy group $\overline{\Phi(G)}_x \cong \cD_x$ is trivial. Thus non-triviality of $\cD_x$ is seen as an obstruction to the transitivity of the action of the automorphism group of $(G_\infty,G)$, and for this reason it was called the \emph{discriminant group} in \cite{DHL2016a}. The article \cite{DHL2016a} also contains examples of actions with finite non-trivial discriminant group, and examples where the discriminant group is a Cantor group.

\medskip

The cardinality of the discriminant group of $(X,G,\Phi)$ may change if one restricts the action to a clopen subset of $X$, see examples in \cite{DHL2016a}. That is, the cardinality of $\cD_x$ reflects the global properties of the action. The idea on how to use the Ellis group to develop a \emph{local} invariant stems from \cite{DHL2016c} and was fully developed in the joint paper with Hurder \cite{HL2017a}. The idea is as follows.

The existence of a homeomorphism $\phi: X \to G_\infty$ allows us to define a series of clopen partitions $\cC_n$ of $X$, where the sets in the $n$-th partition are in bijection with the cosets in $G/G_n$. The action of $G$ on $X$ permutes the sets in $\cC_n$, and the bijection commutes with the action of $G$ on $\cC_n$ and on $G/G_n$. We denote by $X_n$ the set in the partition $\cC_n$ which contains the point $x$. Then the restricted action $\Phi_n = \Phi|_{X_n}$ of $G$ on $X_n$ is the action of the subgroup $G_n$, which is minimal and equicontinuous.

\begin{ex}\label{ex-clopens}
{\rm
Let $T$ be a $d$-ary tree with minimal and equicontinuous action of a discrete group $G$ as in Example \ref{ex-treecylinder}. Let $X = \cP_d$ be the space of infinite paths in $T$.

Let $x = (v_n)_{n \geq 0}$ be a path. Let $G_n = \{g \in G \mid g \cdot v_n = v_n\}$ be the subgroup of elements in $G$ which fix the vertex $v_n$, called the \emph{stabilizer} of $v_n$, or the \emph{isotropy subgroup} of the action of $G$ at $v_n$. Since $G$ acts transitively on the finite set $V_n$, then we have $|G: G_n| = |G/G_n| = |V_n| = d^n$, so $G_n$ has finite index in $G$. If $g \in G$ fixes $v_n$, then it fixes $v_i$ for $ 0 \leq i < n$, which implies that $G_n \subset G_{i}$ for $0 \leq i < n$. So the isotropy subgroups form a nested chain $\{G_n\}_{n \geq 0}$ of finite index subgroups of $G$.

For $n \geq 0$ let $X_n = U_n(v_n) = \{(y_i) \in X \mid y_n = v_n\}$, that is, $X_n$ contains all infinite paths through the vertex $v_n \in V_n$. If $g \in G_n$ and $(y_i) \in X_n$, then $g$ fixes the $n$-th vertex $y_n = v_n$, but may act non-trivially on vertices at the levels $i>n$. Thus the restriction of the action of $G$ to $X_n$ is given by the action of the isotropy subgroup $G_n$. This concludes Example \ref{ex-clopens}.
}
\end{ex}

We now consider the collection of restricted group actions $(X_n,G_n,\Phi_n)$. For each action there is an Ellis group $\overline{\Phi_n(G_n)}$ and the isotropy group $\overline{\Phi_n(G_n)}_x$ of the Ellis group action at $x$. Theorem \ref{thm-groupchains} allows us to obtain representations of the Ellis groups and the isotropy subgroups as inverse limits of finite groups. For each $n \geq 0$ denote the inverse limits by $E_n$ and $\cD^n_x$ respectively. As we explain in more detail in Section \ref{sec-asymptotic}, there are natural maps
    \begin{align}\label{discr}\xymatrix{\cD^0_x \ar[r]^{\psi_{0,1}} & \cD^1_x  \ar[r]^{\psi_{1,2}} & \cD^2_x   \ar[r] & \cdots} \end{align}
between the discriminant groups, which are surjective group homomorphisms. These homomorphisms need not be injective.

\begin{defn}\label{defn-asymptdiscr}
The \emph{asymptotic discriminant} of an equicontinuous minimal action $(X,G,\Phi)$ is the equivalence class of the chain of surjective homomorphisms \eqref{discr}, with respect to the equivalence relation defined in Section \ref{sec-asymptotic}.
\end{defn}

The asymptotic discriminant of an equicontinuous minimal group action was introduced and studied in \cite{HL2017a} and it was shown to be an invariant of return equivalence of group actions; see \cite{CHL2017} for a rigorous definition and properties of return equivalence. Again, even though \cite{HL2017a} makes an assumption of finite generation of the group $G$, the proofs translate verbatim to the case of a countably generated $G$. The introduction of the asymptotic discriminant allowed us to divide equicontinuous minimal Cantor actions into two large classes, as in the following definition.

\begin{defn}\label{defn-stable}
An equicontinuous minimal group action $(X,G,\Phi)$ is \emph{stable} if there exists an integer $m\geq 0$ such that for all $n \geq m$ the group homomorphisms $\psi_{n,n+1}: \cD^n_x \to \cD^{n+1}_x$ in the asymptotic discriminant \eqref{discr} are isomorphisms. If such an $m$ does not exist, then the action $(X,G,\Phi)$ is \emph{wild}.
\end{defn}

If the action is stable (wild), we also say that its asymptotic discriminant is \emph{stable} (\emph{wild}). 

An uncountable family of distinct wild actions of subgroups of $SL(k,\mZ)$ for $k \geq 3$ was constructed in \cite{HL2017a}.  In \cite{DHL2016c}, every finite group and every separable profinite group were realized as discriminant groups of stable actions of subgroups of $SL(k,\mZ)$, $k \geq 3$. In those examples, given a finite or a separable profinite group $D$, we constructed an action of a subgroup of $SL(k,\mZ)$ such that for every $n \geq 0$ in the sequence \eqref{discr} we have $\cD_x^n \cong D$, and every map $\psi_{n,n+1}$ is an isomorphism. 

If the discriminant group $\cD_x^0$ in \eqref{discr} is trivial, then for every $n \geq 0$ the discriminant group $\cD_x^n$ is trivial, and the action $(X,G,\Phi)$ (and all restricted actions $(X_n,G_n,\Phi_n)$) is free.  The converse need not be true, that is, it is possible to construct a free minimal equicontinuous action of a discrete group on a Cantor set which has a non-trivial discriminant group. The first example of an action with this property was obtained in \cite{FO2002}. Stable actions of subgroups of $SL(k,\mathbb{Z})$, described in the previous paragraph are also examples of free actions.

\medskip
The last background concept we need to review before stating our main theorems is arboreal representations of the absolute Galois group of a field $K$.

For a field $K$, let $K^{\rm sep}$ be a separable closure of $K$.  Then the absolute Galois group ${Gal}(K^{\rm sep}/K)$ of $K$ is the group of field automorphisms of $K^{\rm sep}$ which fix $K$. Let $f \in K[x]$ be a polynomial of degree $d \geq 2$. Denote by $f^n = f \circ \cdots \circ f$ the $n$-th iterate of $f$. For a fixed $\alpha \in K$, the polynomial equation $f^n(x) = \alpha $ has $d^n$ zeros, counting multiplicities, in an algebraic closure $\oK$ of $K$. We will assume that for any $n \geq 1$, the $d^n$ zeros are distinct, and so they are contained in the separable closure $K^{\rm sep}$. 

We are now going to build a tree of preimages of $\alpha$ under the iterations $f^n$. To this end, let $V_0 = \{\alpha\}$, let $V_n = f^{-n}(\alpha)$ for $n \geq 1$, that is, $V_n$ is the set of the solutions of the polynomial equation $f^n(x) = \alpha$. Let $V = \bigsqcup_{n \geq 0} V_n$ be the set of vertices of an infinite tree $T$, where we connect $b \in V_{n+1}$ and $c \in V_{n}$ by an edge if and only if $f(b) = c$. Since all zeros of $f^n$ are distinct, every vertex $c \in V_{n}$ is connected to precisely $d$ vertices in $V_{n+1}$, so the tree $T$ is a  $d$-ary rooted tree as in Example \ref{ex-treecylinder}. We denote by $T_n$ the finite subtree of $T$ with vertices $V_0 \cup \cdots \cup V_n$. Denote by $Aut(T)$ and $Aut(T_n)$ the automorphism groups of $T$ and $T_n$ respectively. It is well-known (see, for instance, \cite{BOERT1996}), that $Aut(T_n) = [S_d]^n$, where $[S_d]^n$ denotes the $n$-fold wreath product of symmetric groups $S_d$, and ${\ds Aut(T) = {\lim_{\longleftarrow}} [S_d]^n \cong [S_d]^\infty}$.

For $n \geq 1$, let $K_n = K(f^{-n}(\alpha))$, that is, $K_n$ is a separable extension of $K$ obtained by adjoining the solutions of the polynomial equation $f^{n}(x) = \alpha$. Then $Gal(K^{\rm sep}/K)$ acts on $V_n$ by permuting its elements, and the image of ${ Gal}(K^{\rm sep}/K)$ in $Aut(T_n)$ is a subgroup $H_n$ of the wreath product $[S_d]^n$. We assume that $f^n(x) - \alpha$ is irreducible, so $H_n$ acts transitively on $V_n$. One has $K_{n} \subset K_{n+1}$, and there are natural maps $H_{n+1} \to H_{n}$ \cite[Lemma 4.1]{Odoni1985} , so we obtain a representation 
 \begin{align}\label{arb-1} \rho_{f,\alpha}:  {Gal}(K^{\rm sep}/K) \to H_\infty = \lim_{\longleftarrow}\{H_{n+1} \to H_{n}\} \subseteq {\rm Aut}(T), \end{align}
called the \emph{arboreal representation} of $Gal(K^{sep}/K)$ into the group of automorphisms of the tree $T$. 

\begin{remark}\label{remark-irreducible}
{\rm In the construction above, we made an assumption that every iterate $f^n(x) - \alpha$ is irreducible to ensure that the Galois groups $K(f^{-n}(\alpha))/K$ act transitively on the levels $V_n$ of the tree $T$. Transitivity of the action is required for our main Theorem \ref{thm-suspension}. We note however that an arboreal representation \eqref{arb-1} is defined even if this assumption is not satisfied, although in this case the actions of the Galois groups of the field extensions are considerably more complicated, see the discussion after Theorem 2.4 in \cite{Jones2014}. 
}
\end{remark}

The study of arboreal representations of Galois groups was initiated by Odoni \cite{Odoni1985}, motivated by the relation of arboreal representations to certain questions about the density of primes. One of the questions posed by Odoni was to identify fields and polynomials, for which the representation \eqref{arb-1} is surjective. Recently the study of this and related questions has seen an outburst of activity, see, for example, Jones \cite{Jones2014} for a recent survey. An extensive list of open questions can be found at \cite{AimPL,AIMWorkshop}. 

In this paper, we show that arboreal representations give rise to minimal equicontinuous actions of countably generated discrete groups on Cantor sets. Thus we obtain a new class of examples of such actions. Moreover, we compute the asymptotic discriminant, introduced in \cite{HL2017a} in the topological setting, for these actions.  

A remark on terminology is in order. The reader should remember that the asymptotic discriminant of Cantor group actions, studied here, is a completely different notion to the `discriminant of a polynomial', which is the product of squares of differences of polynomial roots. These two discriminants should not be confused.

\medskip
We are now ready to state our main theorems. 

Recall that $\cP_d$ denotes the space of infinite paths in the $d$-ary tree $T$, which is a Cantor set by Example \ref{ex-treecylinder}.  We denote by ${\bf v} = (v_n)_{n \geq 0}$ a path in $\cP_d$. Given a descending group chain $\{G_n\}_{n \geq 0}$, we denote by ${\bf e}$ the sequence of cosets of the identity in the inverse limit ${\ds G_\infty = \lim_{\longleftarrow}\{G/G_{n+1} \to G/G_n\}}$. 

\begin{thm}\label{thm-suspension}
Let $f(x)$ be a polynomial of degree $d \geq 2$ over a field $K$, and $\alpha \in K$.  
Let $H_\infty$ be the image of the arboreal representation $\rho_{f,\alpha}$ as in \eqref{arb-1}, and let ${\bf v}$ be a path in $\cP_d$. Then there exists a countably generated discrete group $G$, and a group homomorphism
  \begin{align}\label{Phi-actionn} \Phi: G \to Homeo (\cP_d), \end{align}
such that the following is true:
\begin{enumerate}
\item There is an isomorphism 
  \begin{align}\label{eq-isom11}\widetilde{\phi}:\overline{\Phi(G)} \to H_\infty \end{align}
of the Ellis group $\overline{\Phi(G)} \subset Homeo(\cP_d)$  of the action \eqref{Phi-actionn} onto the image of the arboreal representation \eqref{arb-1}.
 
\item There is a descending chain $\{G_n\}_{n \geq 0}$ of subgroups of finite index in $G$, such that there is a homeomorphism 
  \begin{align}\label{eq-isom22}\phi:  G_\infty = \lim_{\longleftarrow}\{G/G_{n+1} \to G/G_{n}\} \to \cP_d\end{align}
of Cantor sets, with $\phi({\bf e}) = {\bf v}$, and for all ${\bf u} \in \cP_d$ and all ${\bf g} \in \overline{\Phi(G)}$ we have
  \begin{align}\label{eq-equivar2}\widetilde{\phi}({\bf g}) \cdot \phi({\bf u}) = \phi({\bf g}(\bf u)).\end{align}
  \end{enumerate}
\end{thm}

Theorem \ref{thm-suspension} associates a group chain $\{G_n\}_{n \geq 0}$ to an arboreal representation $\rho_{f,\alpha}$, and so allows us to compute the asymptotic discriminant of the action of $G$. Our next theorems provide examples of group actions with non-trivial stable or wild asymptotic discriminant which arise from arboreal representations. In the rest of the paper, we sometimes say a `stable/wild arboreal representation', meaning that the asymptotic discriminant of the discrete group action associated to this arboreal representation is stable/wild. 

One of the current questions in the field of arboreal representations is to identify polynomials, such that the image of the arboreal representation \eqref{arb-1} has finite index in $Aut(T)$. This question is a slight generalization of the original question of Odoni \cite{Odoni1985}, about when such a representation is surjective. Odoni conjectured that polynomials with surjective arboreal representations exist in every degree $d$, and found the first specific example of such a polynomial for $d=2$ and $K = \mQ$. Since then, a flurry of results for polynomials of degree $d=2$ have been obtained, we refer to Jones \cite{Jones2014} for an overview. Examples of polynomials satisfying the Odoni conjecture in all prime degrees $d \geq 3$ for $K= \mQ$ were obtained by Looper \cite{Looper2016}.  Juul \cite{Juul2014} studied this question for rational functions over fields of positive characteristic. Examples in even degrees $d \geq 20$  over $\mQ$ were recently obtained by Kadets \cite{Kadets2018}. Benedetto and Juul \cite{BJ2018} found examples in even degrees $d$ over any number field $K$, and in odd degrees when $K$ does not contain $\mQ(\sqrt{d}, \sqrt{d-2})$. Examples in every positive degree were also obtained by Specter \cite{Specter2018}.

Generically one expects arboreal representations to have finite index in $Aut(T)$. A conjecture by Jones states that if $f$ is a quadratic rational function, then the image of the arboreal representation has finite index in $Aut(T)$ unless one of the four specific obstructions is present, see \cite[Section 3]{Jones2014} and specifically Conjecture 3.11. 
Our second result is that if the image in \eqref{arb-1} has finite index in $Aut(T)$, then the asymptotic discriminant of the associated discrete group action is wild. 

\begin{thm}\label{thm-wild}
Let $K$ be a field, let $f(x)$ be a polynomial of degree $d \geq 2$, and let $\alpha \in K$. Let $\rho_{f,\alpha}$ be the associated arboreal representation, given by \eqref{arb-1}, with image $H_\infty$. Let $G$ be a countably generated discrete group acting on the space of paths $\cP_d$ of the $d$-ary rooted tree $T$, given by Theorem \ref{thm-suspension}. If $H_\infty$ has finite index in $Aut(T)$, then the action of $G$ on the space of paths $\cP_d$ is wild.\end{thm}

The proof of Theorem \ref{thm-wild} is geometric, and it does not require the knowledge of specific polynomials or groups in terms of generators and relations. The proof is based on the relationship of the asymptotic discriminant with the \emph{strong quasi-analytic property} of actions. The strong quasi-analytic (or strong quasi-effective) property for pseudogroup actions was introduced in \cite{ALC2009}. For a special case of group actions on Cantor sets, there is an explicit relationship with the algebraic invariant of the asymptotic discriminant, which was studied in \cite{HL2017a}.

Our next examples use representations of Galois groups in terms of generators and relations.

Denote by $\mQ_p$ the field of $p$-adic numbers, where $p$ is an odd prime, and let $K$ be a finite unramified extension of $\mQ_p$ with residual degree $m$. That is, $[K: \mQ_p] = m$, and the residue field of $K$ is the finite field $F_q$, where $q = p^m$. Such an extension is obtained by adjoining to $\mQ_p$ the $(q-1)$-st roots of unity. An unramified extension of $K$ of degree $s >0$ is obtained by adjoining to $K$ $(q^s-1)$-st roots of unity. 

Let $L$ be a finite Galois extension of $K$ with ramification index $e$ and residual degree $s$. Then $[L:K] = es $, and the residual field of $L$ is $F_{q^s}$. We assume that $e$ is coprime to $p$, that is, $L$ is a tamely ramified extension of $K$. Denote by $K^{tame}$ the maximal tamely ramified extension of $K$. Then \cite{PeteClark,Neuk2008} we have an isomorphism
  $$Gal(K^{tame}/K) \cong \widehat{BS(1,q)},$$
where $ \widehat{BS(1,q)}$ denotes the profinite completion of a Baumslag-Solitar group
  $$BS(1,q) = \langle \sigma, \tau \, \mid \, \sigma \tau \sigma^{-1} = \tau^q \rangle.$$

\begin{thm}\label{thm-kummerextame}
Let $p$ and $d$ be distinct odd primes.
Let $K$ be a finite unramified extension of $\mQ_p$ with residue field $F_q$, where $q = p^m$, and consider the arboreal representation $\rho_{f,0}$, where
    $$f(x) = (x+p)^d - p.$$
Then there exists a sequence of integers $\{c_n\}$, and a non-decreasing unbounded sequence of integers $\{s_n\}$, such that the Galois group of the extension $K(f^{-n}(0))/K$ is given by 
  \begin{align}\label{eq-Galoisqp1}H_n = \langle \tau, \sigma  \, \mid \, \sigma \tau \sigma^{-1} = \tau^q, \, \tau^{c_nd^n}=1, \sigma^{s_n}=1 \rangle,\end{align}
 so that $H_n = BS(1,q)/C_n$, where
   \begin{align}\label{eq-normalcorecn} C_n = \langle \tau^{c_nd^n}, \sigma^{s_n}\rangle. \end{align} 
 There is a choice of a path ${\bf v} \in \cP_d$ such that the associated group chain $\{G_n\}_{n\geq 0}$ consists of the groups
  \begin{align}\label{eq-isotropyqp1}G_n = \langle \tau^{c_nd^n}, \sigma \rangle \subset BS(1,q), \, \textrm{where } G_0 = G = BS(1,q). \end{align}
For each $n \geq 0$, the discriminant group $\cD^n_{\bf v}$ is a Cantor group, and the homomorphism $\cD^n_{\bf v} \to \cD^m_{\bf v}$ is an isomorphism for all $m > n $,  so the action is stable.
\end{thm}

The example of the polynomial in Theorem \ref{thm-kummerextame} is standard in Galois theory, and representations for Galois groups \eqref{eq-Galoisqp1} can be found in the literature. The main difficulty of Theorem \ref{thm-kummerextame} lies in computing the groups in \eqref{eq-isotropyqp1}, and the asymptotic discriminant. This is the novel aspect of this theorem.
  
 \medskip 
  
Having introduced a procedure of computing the asymptotic discriminant for arboreal representations of Galois groups, one can pose many questions for further work. The most interesting among them, in the author's opinion, is the following. 

The question when an arboreal representation has finite index in $Aut(T)$ is of interest in arithmetic dynamics due to its applications to the problems about density of primes in number theory. According to Theorem \ref{thm-wild} in this paper,  if an arboreal representation has finite index in $Aut(T)$, then it is wild.  The converse implication is most likely not true, that is, there may exist arboreal representations with wild asymptotic discriminant, whose image has infinite index in $Aut(T)$. The reason we conjecture that to be true is that finite index subgroups of $Aut(T)$ are topologically countably (infinitely) generated, see Remark \ref{remark-notfinite}, while there exist examples of actions of finitely generated groups \cite{DHL2016c,HL2017a}, arising in the topological setting, which are wild. It was shown in \cite{HL2017a} that an action is wild if and only if it is not locally completely strongly quasi-analytic (LCSQA), see Section \ref{sec-solenoids} for a precise definition of this property. Thus the following question is natural.

\begin{quest}
Is it possible to relate the LCSQA property of group actions on Cantor sets,  or the absence of that, with questions in number theory?
\end{quest}  

Finding examples of topologically finitely generated arboreal representations with wild asymptotic discriminant is an interesting question in itself. A possible candidate for such an example may be the one in \cite{BFHJY}. 

Other possible problems are the following. Ingram \cite{Ingram2013} showed that if a polynomial $f(x)$ over a finite extension $K/\mQ_p$ of degree not divisible by $p$ has good reduction, and the point $\alpha$ is not in the Julia set of $f(x)$, then the arboreal representation $\rho_{f,\alpha}$ has finite index in a certain subgroup of $Aut(T)$. We conjecture that the methods of Theorem \ref{thm-kummerextame} may be extended to this class of polynomials to show that the associated action of a discrete group $G$ on the path space $\cP_d$, given by Theorem \ref{thm-suspension}, is stable with infinite discriminant group. Another interesting problem is to compute the asymptotic discriminant of the arboreal representation associated to the polynomial $f(x) = (x+p)^d - p$ over $K= \mQ$, for $\alpha = 0$. For this polynomial, the Galois groups of field extensions $K_n$ are known.

Anderson, Hamblen, Poonen and Walton \cite{AHPW2017} considered arboreal representations of the polynomial $f(x) = x^\ell - c$ over a finite extension $K/\mQ_p$. In this case the behavior of the arboreal representation depends on whether $p$ divides $\ell$, and on the valuation of $c$. It would be interesting to investigate the properties of the asymptotic discriminant for these actions.

Finally, there is an interesting link between arboreal representations and the dynamics of self-similar actions of iterated monodromy groups on trees, as in \cite{Nekr}. Namely, take $K=\mC(t)$ to be the field of rational functions over the complex numbers, let $f(z)$ be a polynomial with complex coefficients, and take $\alpha = t$. In this case the Galois group $H_\infty$ of the extension of $K$, obtained by adjoining the roots of $f^n(z) = t$, for $n \geq 1$, is called the \emph{profinite iterated monodromy group} of $f(z)$. By Proposition 6.4.2 in Nekrashevych \cite{Nekr}, attributed by Nekrashevych to R. Pink, this group is isomorphic to the closure of the action of the discrete \emph{iterated monodromy group} associated to a partial self-covering $f: \cM_1 \to \cM$, where $\cM_1 \subset \cM$ are subsets of $\mC$. Thus in this case the discrete iterated monodromy group is a natural choice for a discrete group in Theorem \ref{thm-suspension}. Iterated monodromy groups are well-studied, see \cite{Nekr} for references and main techniques. In the case when $f(z)$ is \emph{post-critically finite}, that is, the critical points of $f(z)$ have finite forward orbits under $f$, the discrete iterated monodromy group is finitely generated, and so the profinite iterated monodromy group is topologically finitely generated. Some interesting cases for quadratic $f(z)$ were studied by Pink \cite{Pink2013}, who showed that actions of different iterated monodromy groups may have conjugate closures in $Aut(T)$. The work in progress \cite{Lukina2018} studies the asymptotic discriminant of actions of iterated monodromy groups in the case when $f(z)$ is post-critically finite. 
  
{\bf Acknowledgements}. The author thanks Nicole Looper for sharing her knowledge about arboreal representations of Galois groups, and for suggesting local fields as a suitable base field for examples.  Thanks are also due to Steve Hurder for his interest in this project, and for comments on the first draft of this article. The author thanks the anonymous referee for suggesting many useful improvements to the exposition of the paper.

\section{The asymptotic discriminant of an equicontinuous Cantor action}\label{sec-solenoids}

In this section, we recall the necessary background on equicontinuous Cantor actions and the asymptotic discriminant. Main references for this section are works \cite{FO2002,DHL2016a,DHL2016c,HL2017a}.

\subsection{Equicontinuous actions on Cantor sets} 

Let $X$ be a Cantor set, that is, a compact totally disconnected perfect metrizable space. Recall \cite[Section 30]{Willard} that a space $X$ is \emph{perfect} if every point $x \in X$ is an accumulation point of a non-constant sequence of points in $X$. In other words, no point of $X$ is isolated.  

Let $D$ be a metric on $X$, and suppose $\Phi: G \to {Homeo}(X)$ defines an action of a countably generated discrete group $G$ on $X$. The action $(X,G,\Phi)$ is \emph{equicontinuous}, if for any $\epsilon >0$ there exists $\delta>0$ with the following property: for any $g \in G$ and any $x,y \in X$ such that $D(x,y) < \delta$ we have $D(\Phi(g)(x),\Phi(g)(y))< \epsilon$.

The action of a discrete group $G$ on the path space of a $d$-ary rooted tree in Example \ref{ex-treecylinder} is an example of a minimal and equicontinuous action on a Cantor set.

\begin{defn}\label{defn-groupchain}
Let $G$ be a countably generated discrete group. A nested descending sequence $ \{G_n\}_{n \geq 0} = G_0 \supset G_1 \supset G_2 \supset \cdots$ of finite index subgroups of $G$ is called a \emph{group chain}.
\end{defn}

Let $\{G_n\}_{n \geq 0}$ be a group chain as in Definition \ref{defn-groupchain}. Then for every $n \geq 0$ the coset space $G/G_n$ is a finite set. Give $G/G_n$ discrete topology, then by a standard argument the inverse limit ${\ds G_\infty = \lim_{\longleftarrow}\{G/G_{n+1} \to G/G_{n}\}}$ is a Cantor set. There is a natural action of $G$ on $G_\infty$ given by \eqref{eq-leftaction}. This action is transitive on every finite coset space $G/G_n$, so by an argument similar to the one in Example \ref{ex-treecylinder} the action $(G_\infty,G)$ is minimal and equicontinuous. We call $(G_\infty,G)$ the \emph{dynamical system associated to a group chain} $\{G_n\}_{n \geq 0}$. 

It turns out that every minimal equicontinuous group action on a Cantor set is conjugate to a dynamical system, associated to a group chain as in Definition \ref{defn-groupchain}. This is a consequence of the following statement, which can be found in \cite{ClarkHurder2013} and \cite[Appendix]{DHL2016a} for the case when $G$ is finitely generated. The proof carries on to the case of countably generated groups verbatim, so we omit the details.

\begin{prop} \label{prop-chain}
Let $X$ be a Cantor set, and let $G$ be a countably generated group acting on $X$. Suppose that the action of $G$ on $X$ is equicontinuous, and let $x \in X$ be a point. Then there is a descending chain of clopen sets $X=V_0 \supset V_1 \supset \cdots$ with $\bigcap V_n = \{x\}$ such that:
\begin{enumerate}
\item \label{prop-1} For each $n \geq 0$ the collection of translates $\{ \Phi(g)(V_n)\}_{g \in G}$ is a finite partition of $X$ into clopen sets. 
\item \label{prop-2} The collection of elements which preserve $V_n$, that is,
 $$G_n = \{g \in G \mid \Phi(g)(V_n)= V_n\},$$
 is a subgroup of finite index in $G$. 
\end{enumerate}
Thus $G=G_0 \supset G_1 \supset G_2 \supset \cdots$ is a descending chain of subgroups of finite index, moreover, there is a homeomorphism
  $$\phi: X \to G_\infty = \lim_{\longleftarrow}\{G/G_{n+1} \to G/G_{n}\}$$
such that $\phi(g \cdot y) = g \cdot \phi(y)$ for all $y \in X$ and $g \in G$, and such that $\phi(x) = {\bf e}$, where ${\bf e} = (eG_n)$ denotes the sequence of cosets of the identity. 
\end{prop}

Example \ref{ex-clopens} in the Introduction shows how to associate a group chain to the action of a group $G$ on a Cantor set of paths in a $d$-ary tree $T$. 

For a given equicontinuous action $(X,G,\Phi)$, the choice of an associated chain $\{G_n\}_{n \geq 0}$ depends on a choice of a point $x \in X$, and on a choice of clopen sets $V_1 \supset V_2 \supset \cdots$, that is, the choice of the group chain $\{G_n\}_{n \geq 0}$ is not unique. 

\begin{ex}\label{ex-differentbasepoint}
{\rm
Let $(X,G,\Phi)$ be an equicontinuous action, $x \in X$ be a point, and $X=V_0 \supset V_1 \supset V_2 \supset \cdots$ be a collection of clopen sets as in Proposition \ref{prop-chain}. Let $y \in X$ be another point. By Property \eqref{prop-1}, for each $n \geq 0$ the translates $\{ \Phi(g)(V_n)\}_{g \in G}$ form a clopen partition of $X$, so there is a (non-unique) element $g_n \in G$ such that $y \in U_n = \Phi(g_n)(V_n)$. The collection of elements which preserve the clopen set $U_n$ is the conjugate subgroup $g_nG_n g_n^{-1}$, and, associated to the point $y$, there is a group chain $\{G_n' = g_n G_n g_n^{-1}\}_{n \geq 0}$. The inverse limit ${\ds G_\infty' = \lim_{\longleftarrow}\{G/G_{n+1}' \to G/G_n'\}}$ is a Cantor set, and $G$ acts on $G_\infty'$ equicontinuously. The actions $(G_\infty,G)$ and $(G_\infty',G)$ are conjugate, that is, there is a homeomorphism $\psi: G_\infty \to G_\infty'$ such that $\psi(g \cdot y) = g \cdot \psi(y)$ for all $y=(g_nG_n) \in G_\infty$ and all $g \in G$. This isomorphism need not preserve the basepoint, that is, the image $\psi({\bf e})$ of the sequence ${\bf e} = (eG_n)$ of cosets of the identity need not be equal to $(eG_n')$.
}
\end{ex}

Rogers and Tollefson \cite{RT-1} suggested to use the following notion of equivalence of group chains to study the question when two group chains define conjugate actions.

\begin{defn} \cite{RT-1}\label{defn-greq}
In a group $G$, two group chains $\{G_{n}\}_{n \geq 0}$ and $\{H_{n}\}_{n \geq 0}$ with $G_0=H_0=G$ are \emph{equivalent} if there is a group chain $\{K_{n}\}_{n \geq 0}$ and infinite subsequences $\{G_{n_k}\}_{k \geq 0}$ and $\{H_{j_k}\}_{k \geq 0}$ such that $K_{2k} = G_{n_k}$ and $K_{2k+1} = H_{j_k}$ for $k \geq 0$.
\end{defn}

Intuitively, two group chains are equivalent if they can be `intertwined' to form a single descending group chain. Fokkink and Oversteegen \cite{FO2002} investigated equivalence of group chains, determining when two group chains correspond to conjugate actions. A detailed proof of their result, stated below, can also be found in \cite{DHL2016a}.

\begin{thm}\label{thm-equivalence-gc}\cite{FO2002}
Let $(X,G,\Phi)$ and $(X, G, \Psi)$ be equicontinuous actions on a Cantor set $X$ with associated group chains $\{G_n\}_{n \geq 0}$ and $\{H_n\}_{n \geq 0}$, $G_0 = H_0 = G$. Then the actions $(X,G,\Phi)$ and $(X,G,\Psi)$ are conjugate if and only if there exists a sequence of elements $(g_n) \in G$ such that $\{g_n G_n g_n^{-1}\}_{n \geq 0}$ and $\{H_n\}_{n \geq 0}$ are equivalent group chains.
\end{thm}

In Theorem \ref{thm-equivalence-gc}, the elements $g_n$ satisfy the condition $g_nG_m = g_m G_m$ for $n \geq m$, which ensures that the chain $\{g_n G_n g_n^{-1}\}_{n \geq 0}$ is nested. By this theorem, to study an action $(X,G,\Phi)$ in terms of group chains, it is sufficient to consider the chains of conjugate subgroups $\{g_n G_n g_n^{-1}\}_{n \geq 0}$, where $g_nG_m = g_m G_m$ for $n \geq m$.

 \subsection{Ellis group for equicontinuous actions} \label{subsec-ellis}

The  \emph{Ellis (enveloping) semigroup} associated to a continuous group action $\Phi \colon G \times X \to X$ on a topological space $X$   was introduced in the papers \cite{EllisGottschalk1960,Ellis1960}, and is treated in the books   \cite{Auslander1988,Ellis1969,Ellis2014}.   The construction of the Ellis semigroup is abstract, and it can be a   difficult problem to calculate this group exactly. In the case when the action $(X,G,\Phi)$ is equicontinuous, and $X$ is a Cantor set, the situation is a bit simpler, as the Ellis semigroup turns out to be a group and it can be identified with a subgroup of ${Homeo}(X)$. In this section we briefly recall some basic properties   of the Ellis group for a special case of equicontinuous   minimal  systems on Cantor sets.

 Let $X$ be a metric space, and $G$ be a countably generated group acting on $X$ via the homomorphism $\Phi: G \to {Homeo}(X)$. Suppose the action $(X, G, \Phi)$ is equicontinuous. Then the closure  $\overline{\Phi(G)}  \subset {Homeo}(X)$ in the \emph{uniform topology on maps} is identified with the Ellis group of the action. Each element of $\overline{\Phi(G)}$ is the limit of a sequence of maps in $\Phi(G)$, and we use the notation $(g_n)$ to denote a sequence $\{g_n \mid n \geq 0\} \subset G$ such that the sequence $\{\Phi(g_n) \mid n \geq 0\} \subset {Homeo}(X)$ converges in the uniform topology. 

Assume that the action of $G$ on $X$ is minimal, that is, for any $x \in X$ the orbit ${\Phi(G)}(x)$ is dense in $X$.  Then the orbit of the Ellis group $\overline{\Phi(G)}(x) = X$ for any $x \in X$. That is, the group $\overline{\Phi(G)}$ acts transitively on $X$. Denote the isotropy group of the action at $x$ by 
\begin{align}\label{iso-defn2}
\overline{\Phi(G)}_x = \{ (g_n) \in \overline{\Phi(G)} \mid (g_n) \cdot x = x\},
\end{align}
where $(g_n)\cdot x : = (g_n)(x)$, for a homeomorphism $(g_n)$ in $ \overline{\Phi(G)}$. We then have the natural identification $X \cong \overline{\Phi(G)}/\overline{\Phi(G)}_x$ of left $G$-spaces.

Given an  equicontinuous   minimal Cantor system $(X,G,\Phi)$, the  Ellis group     $\overline{\Phi(G)}$  depends only on the image  $\Phi(G) \subset Homeo(X)$, while the isotropy group $\overline{\Phi(G)}_x$ of the action may depend on the point $x \in X$. 
Since the action of $\overline{\Phi(G)}$ is transitive on $X$, given any $y \in X$,   there is an element $(g_n) \in \overline{\Phi(G)}$ such that $(g_n) \cdot x  =y$. It follows that 
  \begin{align}\label{eq-profiniteconj}
  \overline{\Phi(G)}_y = (g_n) \cdot \overline{\Phi(G)}_x \cdot (g_n)^{-1}  \ .
  \end{align}
This tells us that  the \emph{cardinality} of the isotropy group $\overline{\Phi(G)}_x$ is independent of the point $x \in X$, and so  the Ellis group $\overline{\Phi(G)}$ and the cardinality of $\overline{\Phi(G)}_x$ are invariants of $(X,G,\Phi)$.

By Proposition \ref{prop-chain}, given an equicontinuous minimal Cantor action $(X,G,\Phi)$ and a point $x \in X$, there is a group chain $\{G_n\}_{n \geq 0}$ such that $(X,G,\Phi)$ and $(G_\infty,G)$ are conjugate. We now describe the Ellis group of the action $(X,G,\Phi)$ using group chains. For details we refer to \cite{DHL2016a}.

For every $G_n$ consider the \emph{core} of $G_n$, that is, the maximal normal subgroup of $G_n$ given by
\begin{equation}\label{eq-core}
C_n =  {\rm core}_{G}   \, G_n =  \bigcap_{g \in {G}} gG_ng^{-1}  ~ \subseteq ~ G_n \ .
\end{equation}
Since $C_n$ is normal in $G$, the quotient $G/C_n$ is a finite group, and  the collection  $\{C_n \}_{n \geq 0}$ forms a descending chain of normal subgroups of $G$. The inclusions $C_{n+1} \subset C_n$ induce surjective homomorphisms of finite groups, given by
  $$G/C_{n+1} \to G/C_n: gC_{n+1} \mapsto gC_n.$$
The inverse limit space
\begin{align} \label{cinfty-define}  
  C_{\infty}  & =   \lim_{\longleftarrow} \, \left\{  G/C_{n+1} \to G/C_{n}   \right\}  ~ \subset \prod_{n \geq 0} ~ G/C_{n}   
\end{align} 
is a profinite group. Also, since $G_{n+1} \subset G_n$ and $C_{n+1} \subset C_n$, there are well-defined homomorphisms of finite groups $\delta_n: G_{n+1}/C_{n+1} \to G_n/C_n$, and the inverse limit group
  $$\cD_x = \lim_{\longleftarrow}\{G_{n+1}/C_{n+1} \to G_{n}/C_{n}\}$$
is called the \emph{discriminant group} of this action. 

\begin{thm}\cite{DHL2016a}\label{Ellis-Ci}
The profinite group $C_\infty$ is isomorphic to the Ellis group $\overline{\Phi(G)}$ of the action $(X,G,\Phi)$, and the isotropy group $\overline{\Phi(G)}_x$ of the Ellis group action is isomorphic to $\cD_x$.
\end{thm}

By Proposition \ref{prop-chain}, the group chain $\{G_n\}_{n \geq 0}$ depends on the choice of a point $x \in X$, and on the choice of a sequence of clopen sets $X=V_0 \supset V_1 \supset \cdots$ such that the isotropy group of the action of $G$ on $V_n$ is $G_n$. Since the groups $C_n$ are normal, they do not depend on the choice of $x \in X$, but they may depend on the choice of the clopen sets $\{V_n\}_{n \geq 0}$. For any choice of $x$ and $\{V_n\}_{n \geq 0}$, the group $C_\infty$ is isomorphic to the Ellis group $\overline{\Phi(G)}$, so $C_\infty$ is independent of choices up to an isomorphism. One can think of $C_\infty$ as a choice of `coordinates' for the Ellis group $\overline{\Phi(G)}$.

Similarly, the discriminant group $\cD_x$ does not depend on choices up to an isomorphism. We note that, since $\cD_x$ is a closed subgroup of a compact group $C_\infty$, it can either be finite or an infinite profinite group which is topologically a Cantor set. 

\subsection{The asymptotic discriminant}\label{sec-asymptotic}

Let $(X,G,\Phi)$ be a group action on a Cantor set, let $x$ be a point and let $\{G_n\}_{n \geq 0}$ be an associated group chain, that is, the actions $(X,G,\Phi)$ and $(G_\infty,G)$ are conjugate. Recall from Proposition \ref{prop-chain} that the groups $G_m$, $m \geq 0$, preserve the clopen sets $V_m$, that is, for each $m \geq 0$ the restricted action $\Phi_m=\Phi|_{V_m}$ is the action of $G_m$. 

Set $X_m = V_m$, and consider a family of equicontinuous group actions $(X_m,G_m,\Phi_m)$. Then for each $m\geq 0$ we can compute the Ellis group of the action, and the isotropy group at $x$ as follows.
 
For each $n  \geq m \geq 0$, compute the maximal normal subgroup of $G_n$ in $G_m$ by 
 \begin{equation}\label{eq-Mn}
C^m_{n} =   {\rm core}_{G_m}   \, G_{n}  \equiv  \bigcap_{g \in {G_m}} gG_{n} g^{-1}   \subset G_m \ .
\end{equation}

Note that $C^m_{n}$ is the kernel of the action of $G_m$ on the quotient set $G_m/G_{n}$, and $C^0_{n}= C_{n}$.
Moreover,   for all $n > k \geq m \geq 0$, we have $\ds C^m_{n}  \subset C^k_{n} \subset G_{n} \subset G_k \subset G_m$, and $C^m_{n}$ is a normal subgroup of $G_k$. 

\begin{ex}\label{ex-clopens-1}
{\rm
As in Example \ref{ex-clopens}, consider a $d$-ary tree $T$ and let $X = \cP_d$ be the space of paths in $T$. Let $G$ be a discrete group acting minimally and equicontinuously on $X$. Let $(v_n)_{n\geq 0} \in X$ be an infinite path, and denote by $U_n(v_n) = \{(y_i) \in \cP_d \mid y_n = v_n\}$ a clopen subset of paths containing the vertex $v_n$. As in Example \ref{ex-clopens}, let $G_n = \{g \in G \mid g \cdot v_n = v_n\}$ be the isotropy subgroup of the action of $G$ on $V_n$ at $v_n$. Set $X_n = U_n(v_n)$, and $\Phi_n = \Phi|_{X_n}$. Then we have a collection of restricted actions $(X_n,G_n,\Phi_n)$.

For every $n \geq 0$ the maximal normal subgroup $C_n$ of $G_n$ in $G$ consists of elements which fix every vertex in $V_n$. Let $m >0$, and consider $(X_m,G_m,\Phi_m)$. The paths in $X_m = U_m(v_m)$ are contained entirely in a subtree $T_{ \geq v_m}$ of $T$, whose root is the unique path between $v_0$ and $v_m$. The group $G_m$ permutes the paths in $T_{\geq v_m}$ while fixing the root. For $n \geq m$, the maximal normal subgroup $C^m_n$ of $G_n$ in $G_m$ acts trivially on the vertices in the intersection $V_n \cap T_{\geq v_m}$, possibly permuting non-trivially the vertices in $V_n$ which are not in the subtree $T_{\geq v_m}$. Thus every element of $G_n$ which is in $C_n$ is also in $C^m_n$, but elements of $C^m_n$ need not be in $C_n$. So we have $C_n \subseteq C^m_n$.
}
\end{ex}
Define the profinite group 
    \begin{eqnarray}  
C^m_{k,\infty}   & \cong &  \lim_{\longleftarrow} \, \left\{G_{k}/C^m_{n+1} \to G_k/C^m_{n }  \mid n \geq k  \right\} \\ & & = \{ (g_{n} C^m_{n})   \mid  n \geq k   \ , \  g_{k} \in G_{k} \ , \  g_{n +1} C^m_{n} = g_{n} C^m_{n}  \}     \label{Dinfty-definen}  \ .  
   \end{eqnarray}
Then   $C^m_{m,\infty}$ is the Ellis   group of the action $(X_m,G_m,\Phi_m)$, with an associated group chain $\{G_n\}_{n \geq m}$. In particular, $C^0_{0,\infty} = C_{\infty}$, defined by \eqref{cinfty-define} .

Since $G_k \subset G_m$, by definition  we have that $C^m_{k,\infty} \subset C^m_{m,\infty}$, and so  $C^m_{k,\infty}$ is a clopen neighborhood of the identity in $C^m_{m,\infty}$. 

The discriminant group  associated to the truncated group chain   $\{G_n\}_{n \geq m}$ is given by
\begin{eqnarray} 
\cD_x^m & = &    \lim_{\longleftarrow}\, \left \{  G_{n+1}/C^m_{n+1}  \to G_{n}/C^m_{n} \mid n \geq m\right\}    \subset C^m_{m,\infty} \label{eq-discquotients1} \\  & = &    \lim_{\longleftarrow}\, \left \{  G_{n+1}/C^m_{n+1}  \to G_{n}/C^m_{n} \mid n \geq k\right\}    \subset C^m_{k,\infty} , \label{eq-discquotients2}
\end{eqnarray}
where we have $\cD_x^m \subset C^m_{k,\infty}$ since $G_n \subset G_k$ for $n \geq k$. The last statement can be rephrased as saying that the discriminant group $\cD_x^m$ is contained in any clopen neighborhood of the identity in $C^m_{m,\infty}$.

To relate the discriminant groups $\cD_x^m$ and $\cD_x^k$ for $k \geq m$, we define the following maps.

For each $n \geq k \geq m \geq 0$, the inclusion $C^m_{n} \subset C^k_{n}$ induces surjective  group homomorphisms  
\begin{align}\label{eq-quotiensmn} 
\phi_{m,k}^{n} \colon   G_{n}/C^m_{n}  \longrightarrow  G_{n}/C^k_{n} \ ,
\end{align}
and the standard methods show that the maps in \eqref{eq-quotiensmn}  yield   surjective homomorphisms  of the clopen neighborhoods of the identity in $C^m_{m,\infty}$ onto the profinite groups $C^k_{k,\infty}$,
   \begin{align}\label{eq-chomeomorph} 
 \widehat{\phi}_{m,k} \colon   C^m_{k,\infty}  \to C^k_{k,\infty},
 \end{align} 
 which commute with the left action of $G$. Let $\cD_{m,k} \subset C^m_{k,\infty}$ denote the image of $\cD_x^m$ under the map \eqref{eq-chomeomorph}.
 It then follows from  \eqref{eq-quotiensmn}    that for $k > m \geq 0$,  there are    surjective homomorphisms, 
\begin{equation}\label{eq-discmapsnm2}
   \cD_x = \cD_x^0   ~ \stackrel{~  \widehat{\phi}_{0,m} ~ }{\longrightarrow} ~  \cD_{0,m} \cong \cD_x^m ~ \stackrel{~  \widehat{\phi}_{m,k} ~}{\longrightarrow} ~  \cD_x^k \ .
\end{equation}
 
 Thus, given an equicontinuous group action $(X,G,\Phi)$ on a Cantor set $X$, there is an associated sequence of surjective homomorphisms of discriminant groups \eqref{eq-discmapsnm2}, associated to the sequence of truncated group chains $\{G_n\}_{n \geq m}$, $m \geq 0$. 
 
 We now define an equivalence relation on sequences of group homomorphisms as in \eqref{eq-discmapsnm2}, called the \emph{tail equivalence}, first introduced by the author in the joint work with Hurder \cite{HL2017a}.
 
  \begin{defn} \label{def-tailequiv}\cite{HL2017a}
 Let $\cA = \{\phi_{n} \colon A_{n} \to A_{n+1} \mid n \geq 1\}$ and $\cB = \{\psi_{n} \colon B_{n} \to B_{n+1} \mid n \geq 1\}$ be two  sequences of  \emph{surjective}   group homomorphisms. 
 We say that $\cA$ and $\cB$ are tail equivalent, and write $\cA \tail \cB$, if the sequences of groups  $\cA$ and $\cB$ are intertwined 
 by a   sequence of surjective   group homomorphisms. 
 That is, there exists:
 \begin{enumerate}
\item an increasing sequence of indices $\{p_n   \mid n \geq   1 ~ , ~ p_{n+1} > p_{n} \geq n \geq 1 \}$;
\item an increasing sequence of indices $\{q_n   \mid n \geq   1 ~ , ~ q_{n+1} > q_{n} \geq n \geq 1 \}$;
\item a sequence $\cC = \{ \tau_n \colon C_n \to C_{n+1} \mid n \geq 1\}$ of surjective group homomorphisms;
\item a collection of   isomorphisms $\Pi_{\cA\cC} \equiv \{ \Pi^{n}_{\cA \cC} \colon A_{p_n} \to C_{2n-1}\mid n \geq  1\}$;
\item a collection of   isomorphisms $\Pi_{\cB\cC} \equiv \{ \Pi^{n}_{\cB \cC} \colon B_{q_n} \to C_{2n}\mid n \geq  1\}$;
\end{enumerate}
such that  for all $n \geq 1$, the following diagram commutes:
   
    \begin{align} \label{eq-commutativediagram}
 \xymatrix{
 \cdots A_{p_n}   \ar[r]^{\phi_{p_n}} \ar[d]_{\Pi^n_{\cA\cC}}^{\cong} &\cdots  \ar[r]^{\phi_{p_{n+1} -1}} &    A_{p_{n+1}} \ar[r]^{\phi_{p_{n+1}}} \ar[d]_{\Pi^{n+1}_{\cA\cC}}^{\cong} & \cdots   \ar[r]^{\phi_{p_{n+2} -1}}&   A_{p_{n+2}} \ar[r]^{\phi_{p_{n+2}}}  \ar[d]_{\Pi^{n+2}_{\cA\cC}}^{\cong} &     \cdots \\
\cdots C_{2n-1}   \ar[r]_{\tau_{2n-1}}    & C_{2n}   \ar[r]_{\tau_{2n}}  &           C_{2n+1} \ar[r]_{\tau_{2n+1}}   & C_{2n+2}   \ar[r]_{\tau_{2n+2}} & C_{2n+3}    \ar[r]_{\tau_{2n+3}}  &  C_{2n+4}  & \cdots      \\
\cdots     \ar[r]_{\psi_{q_n -1}}    & B_{q_{n }}   \ar[r]_{\psi_{q_n }} \ar[u]_{\Pi^{n}_{\cB\cC}}^{\cong} &     \cdots  \ar[r]_{\psi_{q_{n+1} -1}}   &      B_{q_{n+1}} \ar[r]_{\psi_{q_{n+1}}}    \ar[u]_{\Pi^{n+1}_{\cB\cC}}^{\cong} & \cdots  \ar[r]_{\psi_{q_{n+2} -1}}  & B_{q_{n+2}}  \ar[u]_{\Pi^{n+2}_{\cB\cC}}^{\cong} &  \cdots 
} 
\end{align}
 \end{defn}

A sequence $\cA$ is \emph{constant} if   each map $\phi_{n} \colon A_{n} \to A_{n+1}$ is an isomorphism for all $n \geq 1$, and $\cB$ is said to be \emph{asymptotically constant} if it is tail equivalent to a constant sequence $\cA$. The following  result follows from the usual method of ``chasing of diagrams''.

\begin{lemma}\label{lem-asymptoticconstant}\cite{HL2017a}
A sequence $\ds \cB = \{\psi_{n} \colon B_{n} \to B_{n+1} \mid n \geq 1\}$ of  surjective homomorphisms is asymptotically constant if and only if there exists $n_0 \geq 1$ such that $ker(\psi_{n})$ is trivial for all $n \geq n_0$.
\end{lemma}

We now use the `tail equivalence' of group chains of Definition \ref{def-tailequiv} to introduce the notion of the asymptotic discriminant of a minimal Cantor action, and the notions of a stable and a wild action.

 \begin{defn}\label{def-asympdisc2}\cite{HL2017a}
 Let $(X,G,\Phi)$ be an action of a countably generated group $G$ on a Cantor set $X$, and let $\{G_{n}\}_{n \geq 0}$ be an associated group chain. Then the \emph{asymptotic discriminant} for $\{G_{n}\}_{n \geq 0}$ is the tail equivalence class $[\cD_x^m]_{\infty}$ of   
 the sequence of surjective group homomorphisms 
 \begin{equation}\label{eq-surjectiveDn}
[\cD_x^m]_{\infty} = \{\psi_{m,m+1} \colon \cD_x^m \to \cD_x^{m+1} \mid m \geq 1\}
\end{equation}
 defined by the discriminant groups $\cD_x^m$ for the restricted actions of $G_m$ on the clopen sets $X_m \subset X$. 
 
 The action $(X,G,\Phi)$ is \emph{stable} if the asymptotic discriminant $[\cD_x^m]_{\infty}$ is asymptotically constant, and the action is \emph{wild} otherwise.
  \end{defn}
  
It was shown in \cite{HL2017a} that the asymptotic discriminant is invariant under \emph{return equivalence} of group actions. Intuitively, two actions $(X_1,G,\Phi)$ and $(X_2,G,\Psi)$ are return equivalent, if there are clopen sets $U\subset X_1$ and $V \subset X_2$ such that the collections of local homeomorphisms of $U$ and $V$, induced by the actions, are compatible in a sense made precise in \cite{CHL2017}. This notion is analogous to the notion of Morita equivalence for groupoids. In the joint work with Hurder \cite{HL2017a}, the author constructed an uncountable number of minimal Cantor actions of the same subgroup of ${SL}(n,\mZ)$ with pairwise distinct asymptotic discriminants. These actions are not return equivalent.

If an action $(X,G,\Phi)$ is stable, then its asymptotic discriminant is asymptotically constant, which means that there exists $m_0 \geq 0$ such that for all $k > m \geq m_0$, the group homomorphisms $\widehat{\phi}_{m,k}$ defined in \eqref{eq-discmapsnm2}, are isomorphisms. It was shown in \cite{DHL2016c} that every finite group, and every separable profinite group $D$, can be realized as the discriminant group of a stable minimal Cantor action, that is, in these examples $\cD_x^k \cong D$ for $k \geq m_0$.

\subsection{Strong quasi-analyticity} In this section we briefly describe the relationship between the stability of the action, as defined in Definition \ref{def-asympdisc2}, and geometric properties of the action. For details, we refer the reader to \cite{HL2017a}.

 The  \emph{strong quasi-analyticity} condition for pseudogroup actions   was introduced by  J.~{\'A}lvarez L{\'o}pez and A.~Candel  as the  \emph{quasi--effective property}  in Definition~9.4  in  \cite{ALC2009}, and in revised form as Definition~2.18 in \cite{ALM2016}.
 This condition was motivated by the search for a property equivalent to the quasi-analyticity condition introduced by Haefliger \cite{Haefliger1985}  for    smooth foliations.  We will not reproduce the definition of \cite{ALC2009} here, as we would like to avoid defining pseudogroups and pseudogroup dynamics. Instead we state a version of this definition, adapted to minimal equicontinuous Cantor actions, as in \cite{DHL2016c}.
  
 \begin{defn}\label{defn-lsqa}
 Let $(X,G,\Phi)$ be a minimal equicontinuous group action. Then $(X,G,\Phi)$ is \emph{locally strongly quasi-analytic (LSQA)} if and only if for each $x \in X$ there exists an open set $U \owns x$, such that the following property holds for any $g \in G$: suppose $V \subset U$ is a clopen subset such that the restriction $g|V$ is the identity map. Then $g|U$ is the identity map.
 \end{defn}
 
 Thus, if the action $(X,G,\Phi)$ is LSQA, then there is a subset $U \subset X$ so that we can uniquely extend the action of every element $g \in G$ from a clopen subset $V \subset U$ to $U$. In other words, the restriction of the action of $G$ to $U$ is locally determined. By Proposition \ref{prop-chain}, for any $U$ there is a clopen subset $X_n \subset U$ such that the set of elements $G_n = \{g \in G \mid g\cdot X_n = X_n\}$ is a subgroup of finite index in $G$. Without loss of generality we then can assume $U = X_n$, and the action $(X,G,\Phi)$ is LSQA if and only if there is an $n \geq 0$ such that the restricted action $(X_n,G_n,\Phi_n)$ is locally determined.
 
 \begin{defn}\label{defn-lcsqa}
 Let $(X,G,\Phi)$ be a minimal equicontinuous action of a group $G$ on a Cantor set $X$. We say that $(X,G,\Phi)$ is \emph{locally completely strongly quasi-analytic (LCSQA)} if and only if for each $x \in X$ the action of the Ellis group $\overline{\Phi(G)}$ is LSQA.
 \end{defn}
 
 \begin{remark}
 {\rm
 Let $\Phi:G \to Homeo(X)$ be a group action on a metric space $X$. If $X$ is compact, and the action is equicontinuous, then the uniform topology on maps and the compact-open topology on $Homeo(X)$ coincide. If the image $\Phi(G)$ is compact, then $\Phi(G)$ is equal to its closure in $Homeo(X)$, and then the concepts of LSQA and LCSQA actions coincide. If $\Phi(G)$ is properly contained in its closure in $Homeo(X)$, as is the case for minimal equicontinuous Cantor actions, and the action is LCSQA, then the action is LSQA. \'Alvarez L\'opez and Candel \cite{ALC2009} gave an example of an equicontinuous \emph{non-minimal} action which is LSQA but not LCSQA. It is not known at the moment if there exist examples of minimal actions which are LSQA but not LCSQA.
 }
 \end{remark}
 
 The relationship between the LCSQA property of group actions and their asymptotic discriminant was studied in \cite{HL2017a}. In particular, Proposition 7.4 of \cite{HL2017a} can be rephrased as follows.
 
 \begin{prop}\label{prop-sqa-asymptotic}\cite[Proposition 7.4]{HL2017a}
 Let $(X,G,\Phi)$ be a minimal equicontinuous action of a countably generated group $G$ on a Cantor set $X$, and $x \in X$ be a point. Let $\{G_n\}_{n \geq 0}$ be an associated group chain, and let $[\cD_x^m]_\infty$ be the asymptotic discriminant of the action. Then the asymptotic discriminant $[\cD_x^m]_\infty$ is asymptotically constant if and only if the action $(X,G,\Phi)$ is LCSQA.
 \end{prop}
  
 Using Proposition \ref{prop-sqa-asymptotic}, in some cases it is possible to determine that a minimal equicontinuous action $(X,G,\Phi)$ is wild without computing the asymptotic discriminant, as in the following example.
 
\subsection{Example: Automorphism group of a spherically homogeneous tree}\label{ex-tree}

Let $T$ be a tree, that is, $T$ consists of the set of vertices $V = \cup_{n \geq 0} V_n$, where $V_n$ is finite, and of the set of edges $E$, where an edge $t = [a,b]$ can join two vertices $a$ and $b$ only if $a \in V_{n-1}$ and $b \in V_{n}$, and every vertex $b \in V_{n}$ is joined to precisely one vertex in $V_{n-1}$. We assume that $|V_0| = 1$, that is, $T$ is \emph{rooted}. We also assume that $T$ is \emph{spherically homogeneous}, that is, there is a sequence of positive integers $(\ell_1, \ell_2, \ldots)$ such that every vertex $v \in V_{n-1}$ is joined by an edge to precisely $\ell_{n}$ vertices in $V_{n}$. Note that in this case the following holds: for $n \geq 1$ let $L_n$ be a set with $\ell_n$ elements, then $V_n \cong V_{n-1} \times L_n$.  We use this property of spherically homogeneous trees to compute the automorphism group of $T$ below.  Note that $V_1 \cong V_0 \times L_1 \cong L_1$.
 Also, we assume that $\ell_n \ne 1$ for infinitely many $n \geq 1$. If $\ell_n = d$ for all $n \geq 1$, where $d \geq 2$, then $T$ is a $d$-ary tree as in Example \ref{ex-treecylinder}.
 
 A path in $T$ is an ordered infinite sequence of vertices $(v_n)_{n \geq 0}$ such that for all $n \geq 0$ we have $v_n \in V_n$, and, for $n \geq 1$, there is an edge $[v_{n-1},v_{n}] \in E$ . Denote by $\cP$ the set of all such paths, and put a topology on $\cP$ as in Example \ref{ex-treecylinder}. Then $\cP$ is a Cantor set. 
 
 An automorphism $g$ of a tree $T$ permutes the vertices within each level $V_n$, while preserving the connectedness of the tree. That is, for all $n \geq 1$ the vertices $v_{n-1} \in V_{n-1}$ and $v_{n} \in V_{n}$ are joined by an edge if and only if $g(v_{n-1}) \in V_{n-1}$ and $g(v_{n}) \in V_{n}$ are joined by an edge. Thus $Aut(T)$ acts on infinite paths in the tree. The automorphism group of a spherically homogeneous tree $T$ is computed in \cite{BOERT1996}.  We recall the computation now.

Denote by $S_{\ell_n}$ the symmetric group on $\ell_n$ elements. The group $S_{\ell_n}$ acts by permutations on the set $L_n$ with $\ell_n$ elements, and we denote by $\cdot$ this action. For every $n \geq 1$, fix a bijection $b_n: V_{n-1} \times L_n \to V_n$. 

Denote by $T_n$ the finite subtree of $T$ with vertex set $V_0 \cup \cdots \cup V_n$, and by $Aut(T_n)$ the automorphism group of $T_n$. The tree $T_1$ has $V_1 = L_1$ vertices at level $1$, so
  $$Aut(T_1) \cong S_{\ell_1}.$$
  Now suppose $Aut(T_{n-1})$ is known. Denote by $f: V_{n-1} \to S_{\ell_n}$ a function which assigns to each vertex $v_{n-1} \in V_{n-1}$ a permutation of $L_n$, and let
  $S_{\ell_n}^{V_{n-1}} = \{f: V_{n-1} \to S_{\ell_n}\}$ be the set of all such functions. Then the wreath product
 \begin{eqnarray} \label{wp-formula}S_{\ell_n}^{V_{n-1}}\rtimes Aut(T_{n-1}) \end{eqnarray}
  acts on $V_{n-1} \times L_n$ by
    \begin{eqnarray}\label{wp-action} (f,s)(v_{n-1},w) = (s(v_{n-1}),f(s(v_{n-1})) \cdot w). \end{eqnarray}
 That is, the action \eqref{wp-action} permutes the copies of $L_n$ in the product $V_{n-1} \times L_n$, while permuting elements within each copy of $L_n$ independently. Pushing forward the action to $V_n$ via the bijection $b_n$, we obtain an isomorphism
   \begin{eqnarray}\label{eq-auttn} Aut(T_{n}) \cong S_{\ell_n}^{V_{n-1}}\rtimes Aut(T_{n-1}) \cong S_{\ell_n}^{V_{n-1}} \rtimes  S_{\ell_{n-1}}^{V_{n-2}} \rtimes \cdots \rtimes S_{\ell_1} .\end{eqnarray}
 In the case when $\ell_k = d$ for $1 \leq k \leq n$, we sometimes write $Aut(T_n) \cong [S_d]^n$ for the $n$-fold wreath product in \eqref{eq-auttn}.

 There are natural epimorphisms $Aut(T_n) \to Aut(T_{n-1})$, induced by the projection on the second component in \eqref{wp-formula}. Thus the automorphism group of the tree $T$ is the profinite group
 \begin{align}\label{eq-invlimitwp} Aut(T) = \lim_{\longleftarrow} \{Aut(T_n) \to Aut(T_{n-1}), n\geq 1\} \cong \cdots S_{\ell_n}^{V_{n-1}} \rtimes  S_{\ell_{n-1}}^{V_{n-2}} \rtimes \cdots \rtimes S_{\ell_1}. \end{align}
If $\ell_n =d$ for $n\geq 1$, then we sometimes write $Aut(T) \cong [S_d]^\infty$.
 
 Since $Aut(T)$ is the inverse limit of a collection of finite groups, indexed by natural numbers, then it is topologically countably generated \cite[Proposition 4.1.3]{Wilson}. More precisely, there is a countable collection of elements $G_0=\{g_1,g_2,\cdots\} \subset Aut(T)$, such that the abstract subgroup $G$ generated by $G_0$ is dense in $Aut(T)$ as a set. Denote also by $G$ this subgroup with discrete topology. Since $G$ is identified with a subgroup of $Aut(T)$, and $Aut(T)$ acts on the space of paths $\cP$, we obtain the action $(\cP,G)$ of a countably generated discrete group $G$ on a Cantor set $\cP$. It is straightforward that $Aut(T)$ is isomorphic to the Ellis group of the action $(\cP,G)$.
 
 We claim that this action is wild in the sense of Definition \ref{def-asympdisc2}. Indeed, the following argument shows that the action $(\cP,G)$ is not LCSQA, and then by Proposition \ref{prop-sqa-asymptotic} it is wild.
 
 To show that the action of $Aut(T)$ is not LCSQA, for every open subset $U$ of $\cP$ we exhibit an element which is the identity on a smaller clopen subset $U' \subset U$, but which is not the identity on $U$. 
 
Consider cylinder sets in $\cP$ of the form $U_n(v) = \{ (y_i)_{i \geq 0} \in \cP \mid y_n = v \}$, where $v \in V_n$ is a vertex, and $n \geq 0$. That is, $U_n(v)$ contains all paths in $T$ passing through the vertex $v$. We denote by $T_{\geq v}$ an infinite connected subtree of $T$ with the following property: a path $(y_i)_{i \geq 0}$ is in $U_n(v)$ if and only if $(y_i)_{i \geq 0}$ is a subset of $T_{\geq v}$. Note that since $U_n(v)$ is a clopen subset of a Cantor set $\cP$, it is also a Cantor set. 
 
 Let $U$ be an open set in $\cP$. Then $U$ contains a clopen set $U_k(v)$ for some $v \in V_k$, $k \geq 0$, so without loss of generality we can assume that $U = U_k(v)$. 
 Denote by $Aut(T_{ \geq v})$ the subgroup of $Aut(T)$ which fixes all paths outside of $U_k(v)$. We now compute $Aut(T_{\geq v})$ explicitly.
 
 For $0 \leq n \leq k$ the intersection $V_n \cap T_{\geq v}$ is a single point, so the
 restriction of $Aut(T_{ \geq v})$ to the level $V_n \cap T_{\geq v}$ is the trivial group. The intersection $T_{\geq v} \cap T_{k+1}$ is a finite graph with $\ell_{k+1}$ vertices at level $k+1$, that is, 
    $$T_{\geq v} \cap V_{k+1} \cong L_{k+1}.$$
 It follows that the restriction of $Aut(T_{\geq v})$ to the subgraph $T_{\geq v} \cap T_{k+1}$ of $T_{\geq v}$ is isomorphic to the symmetric group $S_{\ell_{k+1}}$. Following a similar inductive procedure as the one we used to compute $Aut(T)$, we obtain that the restriction of $Aut(T_{\geq v})$ to the finite subgraph $T_{\geq v} \cap T_{k+2}$ is isomorphic to the wreath product
   $$ S_{\ell_{k+2}}^{L_{k+1} } \rtimes  S_{\ell_{k+1}},$$
 and, inductively,
   \begin{align}\label{eq-ksubtreeauto} Aut(T_{\geq v}) \cong \cdots \rtimes S_{\ell_n}^{L_{n-1} \times \cdots \times L_{k+1}} \rtimes  S_{\ell_{n-1}}^{L_{n-2} \times \cdots \times L_{k+1}} \rtimes \cdots \rtimes S_{\ell_{k+1}}.  \end{align}
 If $\ell_n = d$ for all $n \geq k$, then $Aut(T_{\geq v}) \cong [S_d]^\infty$.  
 
 We now show that $Aut(T_{\geq v})$ contains automorphisms which are the identity on a clopen subset of $U_k(v)$, but are not the identity on $U_k(v)$. Let $m >k$, and let $w \in T_{\geq v} \cap V_m$. Then the set $U_m(w)$ of paths containing $w$ is a clopen subset of $U_k(v)$. The complement $U_k(v)-U_m({w})$ is also clopen.  The group $Aut(T_{\geq w})$ is a subgroup of $Aut(T_{\geq v})$, where non-trivial elements are automorphisms which fix the complement $U_k(v)-U_m({w})$ and act non-trivially on $U_m(w)$. Thus for every clopen set $U=U_k(v)$ we produced a collection of elements $Aut(T_{\geq w})$ which are the identity on a clopen subset of $U$ but do not extend to the identity map on $U$. Therefore, the action of $G$ on $\cP$ is not LCSQA.
 
\begin{remark}\label{remark-notfinite}
{\rm
We note that $Aut(T)$ is not topologically finitely generated. Indeed, if $Aut(T)$ is finitely generated, then its abelianization $Aut(T)^{\rm ab}$ must be finitely generated. But $Aut(T)$ is the infinite wreath product of symmetric groups $S_{\ell_n}$, $n \geq 0$. The commutator subgroup of $S_{\ell_n}$ is the alternating group $A_{\ell_n}$, and $S_{\ell_n}/A_{\ell_n} = \mZ/2\mZ$. It follows that $Aut(T)^{\rm ab} \cong \prod_{n = 1}^\infty S_{\ell_n}/A_{\ell_n}$ is isomorphic to the infinite product of copies of $\mZ/2\mZ$, which is not topologically finitely generated. This argument is the same as in \cite{Bond2010} and in \cite[Theorem 3.1]{Jones2014}.
\endproof

}
\end{remark}

\section{Asymptotic discriminant for arboreal representations}\label{sections-provetwotheorems}

In this section, we prove Theorem \ref{thm-suspension} and Theorem \ref{thm-wild}. 

Recall that in the Introduction we defined the arboreal representation
  \begin{align}\label{arb-2} \rho_{f,\alpha}:  {Gal}(K^{\rm sep}/K) \to H_\infty = \lim_{\longleftarrow}\{H_{n+1} \to H_{n}\} \subset {\rm Aut}(T) \end{align}
of a polynomial $f(x)$ of degree $d \geq 2$ over a field $K$.  Here $T$ is a tree with vertex set $V_0  \cup \cdots \cup V_n \cup \cdots$, where $V_0 = \alpha$ and $V_n = f^{-n}(\alpha)$, $n \geq 1$, that is, the vertices in $V_n$ are identified with the solutions of the equation $f^n(x) = \alpha$. The vertices $a \in V_{n+1}$ and $b \in V_{n}$ are joined by an edge if and only if $f(a) = b$. Throughout the paper, we assume that for each $n \geq 1$ the zeros of $f^n(x) = \alpha$ are all distinct. Then $|V_n| = d^n$, and every vertex in $V_{n}$ is joined by edges to precisely $d$ vertices in $V_{n+1}$. That is, the tree $T$ is a $d$-ary tree as in Example \ref{ex-treecylinder}.

The groups $H_n$ in \eqref{arb-2} act on $V_n$ by permutations. The groups $H_n$ are isomorphic to the Galois groups of the field extensions $K(f^{-n}(\alpha))/K$, obtained by adjoining to $K$ the solutions of the equation $f^n(x) = \alpha$. By our assumption that $f^n(x)-\alpha$ is irreducible over $K$ for all $n \geq 1$, we have that $H_n$ acts transitively on $V_n$.
The inverse limit $H_\infty$ is a profinite group isomorphic to the Galois group of the infinite extension of $K$, obtained by adjoining to $K$ the solutions of $f^n(x) = \alpha$ for all $n \geq 1$. This isomorphism is easily seen to be topological (see \cite[Chapter 3]{Wilson} on the topology of profinite Galois groups).

The proof of Theorem \ref{thm-suspension} below is based on two main ideas. Recall from the introduction that we denote by $\cP_d$ the space of paths of a $d$-ary rooted tree $T$.

First note that the image of the arboreal representation $H_\infty$ is a profinite group, while we are looking for a sequence of subgroups of an infinite discrete group $G$. Since $H_\infty$ is the inverse limit of a collection of finite groups, indexed by natural numbers, then it contains a countably generated dense subgroup $G$ \cite[Proposition 4.1.3]{Wilson}. We give $G$ the discrete topology, and consider the action of $G$ on the space of paths $\cP_d$ in a tree $T$. 
Then we identify the Galois groups $H_n$ with the quotient groups $G/C_n$, from where we can determine the groups $C_n$ and $G_n$.

We repeat the statement of Theorem \ref{thm-suspension} for the convenience of the reader. Recall that, given a descending group chain $\{G_n\}_{n \geq 0}$, we denote by ${\bf e}$ the sequence of cosets of the identity in the inverse limit ${\ds G_\infty = \lim_{\longleftarrow}\{G/G_{n+1} \to G/G_n\}}$. 

\begin{thm}\label{thm-suspension-1}
Let $f(x)$ be a polynomial of degree $d \geq 2$ over a field $K$, and $\alpha \in K$. 
Let $H_\infty$ be the image of the arboreal representation $\rho_{f,\alpha}$ as in \eqref{arb-2}, and let ${\bf v}$ be a path in $\cP_d$. Then there exists a countably generated discrete group $G$, and a group homomorphism
  \begin{align}\label{Phi-action} \Phi: G \to Homeo (\cP_d), \end{align}
such that the following is true:
\begin{enumerate}
\item There is an isomorphism 
  \begin{align}\label{eq-isom1}\widetilde{\phi}:\overline{\Phi(G)} \to H_\infty \end{align}
of the Ellis group $\overline{\Phi(G)} \subset Homeo(\cP_d)$  of the action \eqref{Phi-action} onto the image of the arboreal representation \eqref{arb-2}.
 
\item There is a chain $\{G_n\}_{n \geq 0}$ of subgroups of finite index in $G$, such that there is a homeomorphism 
  \begin{align}\label{eq-isom2}\phi:  G_\infty = \lim_{\longleftarrow}\{G/G_{n+1} \to G/G_{n}\} \to \cP_d\end{align}
of Cantor sets, with $\phi({\bf e}) = {\bf v}$, and for all ${\bf u} \in \cP_d$ and all ${\bf g} \in \overline{\Phi(G)}$ we have
  \begin{align}\label{eq-equivar}\widetilde{\phi}({\bf g}) \cdot \phi({\bf u}) = \phi({\bf g}(\bf u)).\end{align}
  \end{enumerate}
\end{thm}
\proof By \cite[Proposition 4.1.3]{Wilson}, $H_\infty$ contains a countable collection $G_0=\{g_1,g_2,\ldots\}$ of elements, such that the closure of $H = \langle G_0 \rangle$ is $H_\infty$. Let $G$ be $H$ with discrete topology. Then $G$ is a countably generated discrete group, which is identified with a dense subgroup $H$ of $H_\infty$.

Denote by $\cdot$ the action of $H_\infty$ on the space of paths $\cP_d$ of the tree $T$. The identification of $G$ with the subgroup $H \subset H_\infty$ induces the homomorphism \eqref{Phi-action},
  that is, it defines an action of $G$ on $\cP_d$. By Example \ref{ex-treecylinder} the action \eqref{Phi-action} is minimal and equicontinuous.   

Let $\varphi_n: H_\infty \to H_n$ be the natural projection onto $H_n$, and denote by $e_n$ the trivial permutation in $H_n$. Then the normal subgroups 
  $$W_n = ker \, \varphi_n = \varphi_n^{-1}(e_n)$$ form a system of open neighborhoods of the identity $(e_n) \in H_\infty$ \cite[Lemma 2.1.1]{RZ}. The subgroups $W_n$ are clopen of finite index \cite[Lemma 2.1.2]{RZ}. For $n \geq 0$, set
  $$\mathcal{C}_n = H \cap W_n,$$
then $\cC_n$ is a normal subgroup of $H$. Since $H$ is dense in $H_\infty$, the restriction $\varphi_n|_H: H \to H_n$ is a surjective homomorphism, with $ker \, \varphi_n|_H = \cC_n$. Thus $\cC_n$ has finite index in $H$. By a standard argument $\{\cC_n\}_{n \geq 0} = \cC_0 \supset \cC_1 \supset \cC_2 \cdots$ is a descending chain of finite index normal subgroups of $H$. 

The group $G$ is identified with $H$. For each $n \geq 0$ let $C_n$ be the subgroup of $G$, identified with $\cC_n$. Then $\{C_n\}_{n \geq 0}$ is a chain of normal subgroups of finite index in $G$, such that 
  \begin{align}\label{all-quotients}H/\cC_n \cong G/C_n \cong H_n. \end{align}
The isomorphisms \eqref{all-quotients} induce the topological group isomorphism
  \begin{align}\label{eq-cinfty} C_\infty = \lim_{\longleftarrow} \{G/C_{n+1} \to G/C_{n}\} \to H_\infty \end{align}
of the inverse limits.  We now associate to the action a group chain $\{G_n\}_{n \geq 0}$, and show that $C_\infty$ is isomorphic to the Ellis group of the action.
  
Since $G$ acts on $\cP_d$ by permutations of paths, there is a natural projection of $\Phi$ to the action  
  $$\Phi_n: G \to Perm(V_n),$$
where $Perm (V_n)$ denotes the group of permutations of the finite set $V_n$. By construction the image $\Phi_n(G)$ in $Perm (V_n)$ is isomorphic to $H_n$, and $ker \, \Phi_n = C_n$, where $ker \, \Phi_n$ is the set of elements acting on $V_n$ trivially.

Recall that we are given a path ${\bf v} = (v_n)_{n \geq 0}$, and for each $n \geq 0$ consider the isotropy subgroup of the action of $G$ at $v_n$, given by
  $$G_n = \{ g \in G \mid \Phi_n(g)(v_n) = v_n\}.$$
Clearly $C_n \subseteq G_n$. Since $H_n$, and thus $G$, acts transitively on $V_n$, we can naturally identify the coset space $G/G_n$ with $V_n$. To this end, define the bijective map $\phi_n: G/G_n \to V_n$ by
  \begin{align}\label{eq-phin} \phi_n(gG_n)  = x_n \in V_n \textrm{ if and only if }  \Phi_n(g)(v_n) = x_n. \end{align}
It follows from the definition \eqref{eq-phin} that for all $g,h \in G$ we have 
 \begin{align}\label{eq-equivn}\phi_n(hgG_n) = \Phi_n(h)(\phi_n(gG_n)) = h \cdot \phi_n(gG_n). \end{align} 
To see that $C_n$ is the maximal normal subgroup of $G_n$, let a subgroup $K \subset G_n$ be normal. Then for any $g \in G$ and any $k \in K$ we have $k g = g k'$ for some $k' \in K$. Acting by $k$ on a coset $g G_n$ we obtain
   $$ k g G_n = g k' G_n = g G_n,$$
since $K \subset G_n$. Thus elements in the normal subgroup $K$ act trivially on $G/G_n$, and so $K \subseteq C_n$. 

Since $\{G_n\}_{n \geq 0}$ is a chain of isotropy groups, and for each $n \geq 0$, $C_n$ is a maximal normal subgroup of $G_n$, it follows by Theorem \ref{thm-groupchains} that the inverse limit group $C_\infty$ in \eqref{eq-cinfty} is topologically isomorphic to the Ellis group $\overline{\Phi(G)}$ of the action \eqref{Phi-action}. Combining this isomorphism with the one in \eqref{eq-cinfty}, we obtain the isomorphism \eqref{eq-isom1}. 

Taking the inverse limit ${\ds \phi = \lim_{\longleftarrow} \phi_n}$, one obtains a homeomorphism \eqref{eq-isom2}, which satisfies \eqref{eq-equivar} for all ${\bf u} \in \cP_d$ and all ${\bf g} \in \Phi(G)$ since every $\phi_n$ has property \eqref{eq-equivn}. This extends to the condition \eqref{eq-equivar} on the closures of $\Phi(G)$ and $H$.
\endproof

We now give our first examples of arboreal representations, which give rise to wild actions of countably generated groups on a Cantor set $\cP_d$, that is, we prove Theorem \ref{thm-wild}. The proof of the theorem is geometric, it uses the non-LCSQA property of wild actions, and it does not require us to know presentations of Galois groups in terms of generators and relations. We restate Theorem \ref{thm-wild} now for the convenience of the reader.

\begin{thm}\label{thm-wild-1}
Let $K$ be a field, let $f(x)$ be a polynomial of degree $d \geq 2$, and let $\alpha \in K$. 
Let $\rho_{f,\alpha}$ be the associated arboreal representation, given by \eqref{arb-2}, with image $H_\infty$. Let $G$ be the countably generated discrete group, given by Theorem \ref{thm-suspension-1}, acting on the space of paths $\cP_d$ of the $d$-ary rooted tree $T$. If $H_\infty$ has finite index in $Aut(T)$, then the action of $G$ on the space of paths $\cP_d$ is wild.
\end{thm}

\proof In the proof, we use notation of Section \ref{ex-tree}. Recall from Section \ref{ex-tree} that if $T$ is a spherically homogeneous tree with $\ell_k = d$, for $k \geq 1$, then its automorphism group is the inverse limit 
  \begin{align}Aut(T) = \lim_{\longleftarrow}\{Aut(T_n) \to Aut(T_{n-1})\},\end{align}
where $Aut(T_n) \cong [S_d]^n$, and $S_d$ is the symmetric group on $d$ elements. 

Denote by $p_n: Aut(T) \to Aut(T_n)$ the projection, and let $e_n$ be the identity in $Aut(T_n)$. Set
 $$W_n = ker \, p_n = p_n^{-1}(e_n),$$
 then $W_n = ker \, p_n$ is a collection of finite index clopen subgroups of $Aut(T)$, which forms a fundamental system of neighborhoods of the identity in $Aut(T)$ \cite[Lemma 2.1.1]{RZ}. 
 
 The group $W_n$ contains permutations of paths in $\cP_d$, which fix every point in the vertex sets $V_k$ for $0 \leq k \leq n$, and which may act non-trivially on the subtrees $T_{\geq v}$ for $v \in V_n$. The action on distinct $T_{\geq v}$ is independent, and $Aut(T_{\geq v}) \cong [S_d]^\infty$ (see the example in Section \ref{ex-tree}), so we obtain
  \begin{align}\label{eq-unformula} W_n \cong \prod_{v \in V_n} Aut(T_{\geq v}) \cong  \underbrace{[S_d]^\infty \times \cdots \times [S_d]^\infty}_{|V_n| \textrm{ times}}.\end{align}

We want to show that if $H_\infty$ has finite index in $Aut(T)$, then it is a clopen subgroup of $Aut(T)$, and so must contain $W_n$ for some $n \geq 0$. Since $Aut(T)$ is topologically infinitely generated (see Remark \ref{remark-notfinite}), $H_\infty$ having finite index in $Aut(T)$ is not sufficient to conclude that it is open. We must show first that $H_\infty$ is closed. 

The image of the arboreal representation $H_\infty$ is the inverse limit of finite groups $H_n \subset Aut(T_n)$, so it is closed in $\prod_n H_n$ \cite[Lemma 1.1.2]{RZ}. Note that the product topology on $\prod_n H_n$ is the relative topology from the product topology on $\prod_n Aut(T_n)$, that is, $W \subset \prod_n H_n$ is open (resp. closed) if and only if there is an open (resp. closed) set $W' \in \prod_n Aut(T_n)$ such that $W = W' \cap \prod_n H_n$. So since $H_\infty$ is closed in $\prod_n H_n$, it is also closed in $\prod_n Aut(T_n)$. Since $Aut(T)$ is a subspace of $\prod_n Aut(T_n)$ and has relative topology from $\prod_n Aut(T_n)$, $H_\infty$ is also closed in $Aut(T)$.

We now have that $H_\infty$ is a closed subgroup of finite index in $Aut(T)$. Then by \cite[Proposition 2.1.2]{RZ} it is also an open subgroup of $Aut(T)$. A clopen subgroup $H_\infty$ must contain a neighborhood of the identity $W_n$ for some $n \geq 0$. Choose $v \in V_n$, and let $T_{\geq v}$ be a subtree of $T$ containing paths through $v$. The path space of $T_{\geq v}$ is a clopen subset $U_n(v) = \{(y_i) \in \cP_d \mid y_n = v\}$. By \eqref{eq-unformula} $H_\infty$ contains the subgroup $Aut(T_{\geq v})$, which permutes the paths in $U_n(v)$ and fixes $\cP_d - U_n(v)$. By the example in Section \ref{ex-tree} the action of $G \cap Aut(T_{\geq v})$ on $U_n(v)$ is not LCSQA and is, therefore, wild. Thus the action of $G$ on $T$ is wild.

\endproof

\section{Stable arboreal representations}\label{sec-arboreal}

In this section we prove Theorem \ref{thm-kummerextame}. Theorem \ref{thm-kummerextame} computes the asymptotic discriminant for the arboreal representation of the polynomial $f(x) = (x+p)^d - p$ over unramified finite extensions of the $p$-adic numbers $\mQ_p$. So first we briefly recall a few facts about the structure of $p$-adic fields, and about the Galois groups of their finite extensions.

\subsection{Galois groups over finite extensions of $p$-adic numbers}\label{sec-extensions} For the background on the $p$-adic numbers and their finite extensions, described here, we refer to \cite{Serre1973,Shatz1972,Lang1965,Cassels1986,PeteClark}.

Let $p$ be a prime, and ${\displaystyle \mZ_p = \lim_{\longleftarrow}\{\mZ/p^n\mZ \to \mZ/p^{n-1}\mZ\}}$ be the ring of $p$-adic integers. Denote by $\mathbb{U}$ the group of invertible elements in $\mZ_p$, then $x \in \mathbb{U}$ if and only if $x$ is not divisible by $p$. Moreover, for every $x \in \mZ_p$ there is a unique representation $x = p^n u$, with $u \in \mathbb{U}$ and $n \geq 0$. Non-invertible elements form a unique maximal ideal $\mathfrak{m}_p$ in $\mZ_p$ \cite{Serre1973}. A ring with unique maximal ideal is called a \emph{local ring}, so $\mZ_p$ is a local ring. The field $F_p = \mZ_p/\mathfrak{m}_p$ is a field with $p$ elements, called the \emph{residue field} of $\mZ_p$.

The $p$-adic rationals $\mQ_p$ is the field of fractions (the quotient field) of the ring $\mZ_p$. It can be obtained by adjoining to $\mZ_p$ a single element, $\mQ_p = \mZ_p[p^{-1}]$. The rationals $\mQ$ embed into $\mQ_p$, so $\mQ_p$ is a field of characteristic zero. The residue field $F_p = \mZ_p/\mathfrak{m}_p$ has characteristic $p$.
  
Let $K$ be a finite extension of $\mQ_p$. Then $K$ can be obtained by adjoining to $\mathbb{Q}_p$ the roots of a polynomial $p(x)$ in $\mathbb{Q}_p[x]$ \cite[Chapter VII]{Lang1965}. The ring of integers $\cO_K$ in $K$ is defined to be the integral closure of $\mZ_p$ in $K$ \cite[Chapter IX]{Lang1965}.  It turns out that $\cO_K$ contains a unique maximal ideal $\mathfrak{m}_K$, and the residue field $k = \cO_K/\mathfrak{m}_K$ is a field of characteristic $p$ with $q = p^m$ elements, where the exponent $m$ is called the \emph{residual degree} of $K/\mQ_p$. We have $[k: F_p] = m$. An extension $K/\mQ_p$ is \emph{unramified} if $[K:\mQ_p] = [k : F_p]$.

Assume that $K/\mQ_p$ is unramified. The finite field $k$ can be obtained from $F_p$ by adjoining the roots of the polynomial $x^{q} - x$, so $Gal(k/F_p)$ is cyclic, generated by the Frobenius automorphism $\sigma(\zeta) = \zeta^p$, where $\zeta$ is a root of unity in the multiplicative group $k^*$. By Hensel's lemma \cite{Shatz1972,Cassels1986} every $\zeta \in k^*$ lifts to a root of unity in $K$. Moreover, one can show that the group of roots of unity in $K$ of order coprime to $p$ is isomorphic to the multiplicative group of the residue field $k$.  It follows that $K$ is obtained by joining to $\mQ_p$ the $(q-1)$-st roots of unity. The Frobenius automorphism $\sigma$ lifts from $k^*$ to an automorphism of $K/\mQ_p$, also denoted by $\sigma$ \cite{PeteClark}. 

Let $K$ be an unramified extension of $\mQ_p$ of degree $m$, and let $L$ be a finite extension of $K$. As above, $k$ is the residue field of $K$, then $k \cong F_q$, $q = p^m$.  Let $\ell$ be the residue field of $L$, and suppose $[\ell : k] = s$.  Then $\ell$ is a field with $q^s$ elements, obtained by joining to $k$ the $(q^s - 1)$-st roots of unity. 

Let $Gal(L/K)$ be the Galois group of the field extension $L/K$, and let $Gal(\ell/k)$ be the Galois group of the residue field extension. 
One can show that there is a surjective homomorphism
  $$Gal(L/K) \to Gal(\ell/k),$$
whose kernel $I$ is called the \emph{inertia subgroup} of $Gal(L/K)$. The cardinality $e$ of $I$ is called the \emph{ramification index} of the extension $L/K$, and $[L:K] = es$.

An extension $L/K$ is \emph{unramified} if $e = 1$ and so $[L:K] = s$. An extension $L/K$ is \emph{totally ramified} if $s=1$ and so $[L:K] = e$. The extension $L/K$ is \emph{tamely ramified} if the ramification index $e$ is coprime to the residue characteristic $p$.

So suppose the extension $L/K$ under consideration is tamely ramified. Then $L/K$ factors through a unique unramified extension $K_s/K$ with $[K_s : K] = [\ell:k]$. Thus $K_s$ is obtained by adding to $K$ the $(q^s-1)$-st roots of unity.  The Galois group $Gal(K_s/K)$ is cyclic, it is generated by the Frobenius automorphism $\sigma: \zeta \mapsto \zeta^q$, where $\zeta$ is a $(q^s-1)$-st root of unity in $K$.

The extension $L/K_s$ is totally ramified, it is obtained by adding to $K_s$ the roots of the polynomial $x^{q^s-1} - p$. This polynomial satisfies the Eisenstein criterion \cite{Shatz1972,Cassels1986}, so it is irreducible. This implies that the Galois group $Gal(L/K_s)$ is cyclic of order $q^s-1$. More precisely, if $\alpha$ is a root of $x^{q^s-1} - p$, and $\zeta$ is a primitive $(q^s-1)$-st root of unity in $K_n$, then every root of $x^{q^s-1} - p$ in $L$ is of the form $\alpha \zeta^i$, $1 \leq i \leq q^s-1$. A generator $\tau$ of $L/K_n$ is given by $\tau: \alpha \zeta^i \mapsto \alpha\zeta^{i+1}$.

Now consider the Galois group of the extension $L/K$. 
The Frobenius automorphism $\sigma$ lifts from $K_s$ to an automorphism of $L$, given by $\sigma: \alpha \mapsto \alpha, \zeta \mapsto \zeta^q$.  A generator of $Gal(L/K_s)$ is given by $\tau:  \zeta \mapsto \zeta, \alpha \mapsto \alpha \zeta$. Combining that together one obtains that \cite{PeteClark}
  \begin{align}\label{eq-trext} Gal(L/K) = \langle \sigma, \tau \mid \sigma \tau \sigma^{-1} = \tau^q, \, \sigma^s = 1, \, \tau^{q^{s}-1} = 1 \rangle. \end{align}
The maximal unramified extension $K^{unr}/K$ is obtained by adjoining to $K$ the $(q^s-1)$-st roots of unity for all $s \geq 1$, and the maximal tamely ramified extension $K^{tame}/K^{unr}$ by adjoining to $K^{unr}$ the $(q^s-1)$-st roots of $p$, for $s \geq 1$. The extension $K^{tame}/K$ contains all $n$-th roots of $p$ for $n$ coprime to $p$, and there is an isomorphism \cite{PeteClark}
  \begin{align}\label{eq-tamegalois} Gal(K^{tame}/K)  = \widehat{BS(1,q)}, \end{align}  
where  $\widehat{BS(1,q)}$ denotes the profinite completion of the Baumslag-Solitar group
  \begin{align}\label{eq-bsgroup}BS(1,q) = \langle \sigma, \tau \mid \sigma \tau \sigma^{-1} = \tau^q \rangle.\end{align}

\subsection{Arboreal representations over local fields}\label{sec-example}

We now prove Theorem  \ref{thm-kummerextame}.
For the convenience of the reader, we restate the theorem before proving it.

\begin{thm}\label{thm-kummerex-1}
Let $p$ and $d$ be distinct odd primes.
Let $K$ be a finite unramified extension of $\mQ_p$ with residue field $F_q$, where $q = p^m$, and consider the arboreal representation $\rho_{f,0}$, where
    $$f(x) = (x+p)^d - p.$$
Then there exists a sequence of integers $\{c_n\}$, and a non-decreasing unbounded sequence of integers $\{s_n\}$, such that the Galois group of the extension $K(f^{-n}(0))/K$ is given by 
  \begin{align}\label{eq-Galoisqp}H_n = \langle \tau, \sigma  \, \mid \, \sigma \tau \sigma^{-1} = \tau^q, \, \tau^{c_nd^n}=1, \sigma^{s_n}=1 \rangle,\end{align}
 so that $H_n = BS(1,q)/C_n$, where
   \begin{align}\label{eq-normalcorecn} C_n = \langle \tau^{c_nd^n}, \sigma^{s_n}\rangle. \end{align} 
 Also, there is a choice of a path ${\bf v} \in \cP_d$ such that the associated group chain $\{G_n\}_{n\geq 0}$ consists of the groups
  \begin{align}\label{eq-isotropyqp}G_n = \langle \tau^{c_nd^n}, \sigma \rangle \subset BS(1,q), \, \textrm{where } G_0 = G = BS(1,q). \end{align}
For each $n \geq 0$, the discriminant group $\cD^n_{\bf v}$ is a Cantor group, and the homomorphism $\cD^n_{\bf v} \to \cD^m_{\bf v}$ is an isomorphism for all $m > n $,  so the action is stable.
\end{thm}

The presentations of Galois groups \eqref{eq-Galoisqp} can be deduced from the literature, see Section \ref{sec-extensions}. The original contribution of the theorem consists of finding explicit representations of the isotropy groups $G_n$ in \eqref{eq-isotropyqp}, and of computing the discriminant groups $\cD^n_{\bf v}$. One of the components, needed to compute the discriminant groups, is the behavior of the exponents $\{s_n\}$ in \eqref{eq-normalcorecn}. We explain in detail how $s_n$ arise in the proof below.
 
\proof To compute the Galois groups of the iterates $f^n(x)$ of our polynomial we will use the fact that $f(x)$ is conjugate by a linear map to the power map $x^d $. This is well-known in the field of the arboreal representations (see, for example, \cite{Odoni1985}). We reproduce the computation here for readers not familiar with arboreal literature. 

Let $y$ be a solution of $f(x) = 0$, and consider $x$ such that $(x+p)^d - p = y$. Then
  $$ (y+p)^d - p = ( (x+p)^d-p+p)^d - p = (x+p)^{d^2} - p =0,$$
and inductively one obtains that the $n$-th iteration of $f^n(x)$ is given by
  $$f^n(x) = (x + p)^{d^n} - p.$$
Let $y_n$ be a solution of $f^n(x) = 0$, let $z_n = y_n + p$, and note that 
  $$(z_n)^{d^n} = (y_n+p)^{d^n} = p.$$
 It follows that $y_n$ is a solution of $f^n(x) = 0$ if and only if $z_n = y_n + p$ is a solution of $p_n(x) = x^{d^n} - p =0$. Since $p$ is in $K$, adjoining $z_n$ to our ground field $K$ generates the same extension as adjoining $y_n$ to $K$. Thus 
   $$Gal(K(f^{-n}(0))) = Gal(K(p_n^{-1}(0))).$$
  
We now can start the proof of Theorem \ref{thm-kummerex-1}. 
Denote $P_n = K(p_n^{-1}(0))$, that is, $P_n$ is obtained by adjoining to $K$ the $d^n$-th roots of $p$. We are going to find a tamely ramified extension $L_n$ of $K$ that contains $P_n$. 

Recall from Section \ref{sec-extensions} that adding to $K$ the $(q^s-1)$-st roots of $p$ produces a tamely ramified extension of $K$. To find an extension that contains $P_n$, consider the equation 
  \begin{align}\label{eq-divisble} q^x = 1 \mod d^n. \end{align}
That is, we are looking for a value of $x$ such that $q^x-1$ is an integer multiple of $d^n$.  
  
Since $q$ and $d$ are coprime, a solution to \eqref{eq-divisble} always exists, for example, one can take 
  $$x= \phi(d^n) = d^{n-1}(d-1),$$ 
 where $\phi$ is the Euler totient function. Let $s_n>0$ be the smallest number such that \eqref{eq-divisble} holds, that is, $0 < s_n \leq \phi(d^n)$.
Let $c_n >0$ be so that 
  \begin{align}\label{divisibility-cond}q^{s_n} -1 = c_n d^n. \end{align}
Since $s_n$ is the smallest positive solution to \eqref{eq-divisble}, it is immediate that the sequence $\{s_n\}$ is non-decreasing, and  $s_n \to \infty$ as $n \to \infty$.
 
 Let $L_n$ be a tamely ramified extension of $K$ of residual degree $s_n$. Then the Galois group $Gal(L_n/K)$ is given by \eqref{eq-trext} with $s = s_n$. Since $d^n$ divides $q^{s_n}-1$, the $d^n$-th roots of $p$ are contained in $L_n$, so $P_n \subset L_n$. In \eqref{eq-trext}, the generator $\tau$ has order $(q^{s_n}-1)$ and so acts transitively on the $(q^{s_n}-1)$-st roots of $p$. Then the power $\tau^{c_n}$ has order $d^n$, and it acts on the $d^n$-th roots of $p$ by permutations. 
 
 Next, let $\zeta$ be a primitive $(q^{s_n}-1)$-st root of unity in $L_n$.  If $(\zeta^i)^{d^n} = 1$, then $((\zeta^{i})^q)^{d^n} = 1$, so the Frobenius automorphism $\sigma$ permutes the $d^n$-th roots of $p$ in $L_n$, and preserves $P_n$. 
 
 The Galois group $Gal(L_n/K)$ injects into the multiplicative group $(\mZ/d^n\mZ)^*$, onto a subgroup generated by $q$. So the order of $\sigma$ is equal to the order of $q$ in $(\mZ/d^n\mZ)^*$. The order of $\sigma$ is the smallest solution of \eqref{eq-divisble}, and is therefore equal to $s_n$. Thus we obtain \eqref{eq-Galoisqp}.

Now let $H = BS(1,q)$, where $BS(1,q)$ is given by \eqref{eq-bsgroup}. Then $H_n = H/C_n$, where
   \begin{align}\label{eq-Cn} C_n = \langle \tau^{c_nd^n}, \sigma^{s_n} \rangle. \end{align}
Using \eqref{divisibility-cond} and Lemma \ref{normalcore-appendix} in Appendix \ref{sec-normalcore-appendix}, we conclude that $C_n$ is a normal subgroup of $H$. We now choose the path ${\bf v}$ in $\cP_d$ and compute the group chain $\{G_n\}_{n \geq 0}$.

Recall that by construction every vertex in $V_n$ corresponds to a $d^n$-th root of $x^{d^n}-p$, which in its turn corresponds to a $d^n$-th root of unity in the unramified extension $K_n$. The multiplicative identity in $K_n \cap \mQ$ is a $d^n$-th root of unity; 
 let $v_n$ be the vertex in $V_n$ corresponding to this root. Clearly, we have $f(v_{n+1}) = v_n$. Since the Frobenius automorphism $\sigma$ fixes the multiplicative identity,   it follows that the isotropy group of the action of $H$ at $v_n$ is
  \begin{align}\label{eq-isotropy}G_n = \langle \tau^{c_nd^n},\sigma \rangle.\end{align}

We now compute the asymptotic discriminant for the group chain $\{G_n\}_{n \geq 0}$. For $n>m$, denote by $C^m_n$ the maximal normal subgroup of $G_n$ in $G_m$, that is,
   $$C^m_n = \bigcap_{g \in G_m} gG_ng^{-1}.$$

\begin{prop}\label{prop-discrx2} 
For all $m \geq 0$ the discriminant group 
  $${\ds \cD^m_{\bf v} = \lim_{\longleftarrow} \{G_{n+1}/C^m_{n+1} \to G_{n}/C^m_{n} \, \mid \, n \geq m\}}$$ 
is infinite and the homomorphisms  $\cD^m_{\bf v}\to \cD^k_{\bf v}$, $k>m \geq 0$, are isomorphisms. 
\end{prop}

That is, Proposition \ref{prop-discrx2} proves that the action of $H$ on the space of paths $\cP_d$ in the tree $T$ is stable. 

\proof Consider the subgroups $C^m_n$ of $G_n$. By \cite{Dudkin2010}, a finite index subgroup in $BS(1,q)$ can be written down as 
  $$C_{n}^m = \langle \tau^{\ell_{m,n}}, \sigma^{k_{m,n}} \tau^{r_{m,n}}\rangle,$$ 
for some values $r_{m,n},k_{m,n}$ and $\ell_{m,n}$, where $0 \leq r_{m,n} \leq \ell_{m,n}-1$, and $\ell_{m,n}$ is coprime to $q$. That is, $C_n^m$ is generated by two elements, which we now compute. 

First, since $C_n \subset C_{n}^m \subset G_n$, and both $C_n$ and $G_n$ have $\tau^{c_n d^n}$ as a generator, we have that $\ell_{m,n} = c_n d^n$. 

Second, since $C^m_n \subset G_n$, a generator $\sigma^{k_{m,n}} \tau^{r_{m,n}}$ must be representable as the composition of the generators of $G_n$. The generators of $G_n$ are $\tau^{c_nd^n}$ and $\sigma$.  Formulas \eqref{id-2} and \eqref{id-4} in the Appendix show that applying relations to a composition of generators of $G_n$ preserves the number of occurrences of $\sigma$ in the composition, and so imply that $r_{m,n}$ must be a power of $c_nd^n$. Since $0 \leq r_{m,n} \leq c_nd^n-1$, we conclude that $r_{m,n} = 0$. We now have
  \begin{align}\label{eq-cnm}C_{n}^m = \langle \tau^{c_n d^n}, \sigma^{k_{m,n}} \rangle.\end{align}

\begin{lemma}\label{lemma-kgrows}
In \eqref{eq-cnm}, $k_{m,n}$ is non-decreasing and $k_{m,n} \to \infty$ as $n \to \infty$. 
 \end{lemma} 
 \proof
Condition \eqref{st1} in Lemma \ref{normalcore-appendix}  implies that for $C^m_n$ to be a normal subgroup of $G_m$ we must have
  \begin{align}\label{eq-div}q^{k_{m,n}} = 1 \mod \frac{c_n d^n}{c_m d^m}.\end{align}
  Different $k_{m,n}$, satisfying \eqref{eq-div} correspond to different subgroups of $G_n$, and the smallest such $k_{m,n}$ corresponds to the maximal normal subgroup of $G_n$ in $G_m$. 
  
Since the sequence $s_n$ in \eqref{eq-divisble} tends to infinity with $n$, for $n$ large enough  $ \frac{c_n d^n}{c_m d^m} >1$, and the ratio $ \frac{c_n d^n}{c_m d^m}$ is non-decreasing and tends to infinity with $n$. It follows that $k_{m,n} \to \infty$ as $n \to \infty$, and $k_{m,n}$ is non-decreasing.
 \endproof
 
 Note that $k_{0,n} = s_n$, since $C^0_n = C_n$.
 
We are now in a position to compute the discriminant group $\cD^m_{\bf v}$.
We note that the cosets in $G_n/C^m_{n}$ are given by
  $$C^m_{n}, \sigma C^m_{n}, \ldots, \sigma^{k_{m,n} -1}C^m_{n},$$
and the inclusions 
  $$G_{n+1}/C^m_{n+1} \to G_{n}/C^m_{n}: \sigma^\alpha C^m_{n+1} \mapsto \sigma^\alpha C^m_{n}$$ 
are clearly surjective. Since $k_{m,n} \to \infty$ as $n \to \infty$, the discriminant group
  $$\cD^m_{\bf v} = \lim_{\longleftarrow}\{G_{n+1}/C^m_{n+1} \to G_{n}/C^m_{n}\}$$
is an infinite profinite group.

We now have to show that the asymptotic discriminant is constant, that is, the surjective homomorphisms $\cD^m_{\bf v} \to \cD^k_{\bf v}$, $k>m \geq 0$ are isomorphisms. For that it is enough to show that the homomorphisms $\cD^0_{\bf v} \to \cD^m_{\bf v}$ are isomorphisms.

So consider the surjective mappings, induced by coset inclusions $C_n \subset C_n^m$, and given by
  $$h_n: G_n/C_n \to G_n/C^m_{n}: \sigma^\alpha C_n \mapsto \sigma^\alpha C_{n}^m.$$
In the limit, there is a surjective homomorphism
  \begin{align}\label{eq-surjective}h_\infty: \cD^0_{\bf v} \to \cD^m_{\bf v}\end{align}
which we now show to be injective.  

Suppose $(\sigma^{a_n} C_n) \ne (\sigma^{t_n} C_n) \in \cD^0_{\bf v}$. Then there is $r_0> 0$ such that for all $n \geq r_0$ we have $\sigma^{a_n} C_n \ne \sigma^{t_n} C_n$, in other words,  $a_n \ne t_n \mod k_{0,n}$, and in particular $a_n \ne t_n \mod k_{0,r_0}$. 

By Lemma \ref{lemma-kgrows} the sequence $k_{m,n} \to \infty$ as $n \to \infty$ and $k_{m,n}$ is non-decreasing, so we can find a number $N> 0$  such that for all $n>N$ we have $k_{m,n} > k_{0,r_0}$. This implies that 
  $$a_n \ne t_n \mod k_{m,n},$$ 
  and the cosets $\sigma^{a_n} C^m_n $ and $\sigma^{t_n} C^m_n$ are distinct. Thus
   $$h_n(\sigma^{a_n} C_n) \ne h_n(\sigma^{t_n} C_n).$$
Thus \eqref{eq-surjective} is an isomorphism, and the action under consideration is stable. 
\endproof

\appendix

\section{Subgroups of finite index in the Baumslag-Solitar group}\label{sec-normalcore-appendix}

Consider Baumslag-Solitar groups with a presentation
 $$BS(1,q) = \langle \tau, \sigma \, \mid \, \sigma \tau \sigma^{-1} = \tau^q \rangle. $$
By \cite{Dudkin2010}, every subgroup of finite index in $BS(1,q)$  can be represented using generators as
  \begin{align}\label{eq-subgroups} H = \langle \tau^\ell, \sigma^m \tau^s \rangle, \end{align}
where $\ell$ and $q$ are coprime, $1 \leq s \leq  \ell-1$, and $m \ell = k$ is the index of $H$. All subgroups are distinct for different values of $\ell, s$ and $m$.

Also, $H$ is normal in $BS(1,q)$ if and only if these two conditions are satisfied at the same time:
\begin{eqnarray}
 \label{cond1} {\rm a)}\, \ell \, | \, q^m - 1 & \rm{and}& {\rm b)}\, \label{cond2}\ell \, \vert \, s(q-1).
\end{eqnarray}

In this section we make a computation similar to the ones in \cite[Section 1]{Dudkin2010} in order to determine when a finite index subgroup $S$ of $BS(1,q)$ is normal in a proper subgroup $H$ of $BS(1,q)$. The divisibility conditions of \cite{Dudkin2010} can be obtained from ours by setting $H = BS(1,q)$, up to a sign.

\begin{lemma}\label{normalcore-appendix}
Let $t$ and $r$ be coprime to $q$. Let $H = \langle \tau^r, \sigma^{\alpha} \rangle$ be a subgroup of $BS(1,q)$, and let $S = \langle \tau^t, \sigma^m \tau^s \rangle$ be a subgroup of $H$, $0 \leq s \leq t-1$. Then
  \begin{eqnarray}\label{st1} \tau^{r} (\sigma^m \tau^s) \tau^{-r} \in S & \textrm{ if and only if } & \frac{t}{r} \, \textrm{ divides } \, 1- q^m, \\ \label{st2}
    \sigma^{\alpha}(\sigma^m \tau^s)\sigma^{-\alpha} \in S & \textrm{ if and only if } & t \, \textrm{ divides } \, s(q^\alpha-1), \  \\ \label{st3} \sigma^{\alpha} \tau^{t} \sigma^{-\alpha} \in S. & & \end{eqnarray}
  That is, $S$ is normal in $H$ if and only if ${\displaystyle \frac{t}{r}}$ divides $1 - q^m$, and $t$ divides $s(q^\alpha - 1)$.  
\end{lemma}

\proof The relation $ \sigma \tau  = \tau^q \sigma$ generates the identity
\begin{eqnarray} \label{id-2} \sigma^u \tau &= \tau^{q^u} \sigma^u. \end{eqnarray}
Applying \eqref{id-2} $\beta$ times yields the identity
  \begin{align}\label{id-4} \sigma^u \tau^\beta = \tau^{\beta q^u} \sigma^u. \end{align}

Note that $r$ divides $t$, since $S \subset H$, and consider \eqref{st1}. Using \eqref{id-4} we obtain
\begin{align*}\label{eq-11}\tau^{r} (\sigma^m \tau^s) \tau^{-r}= \tau^{r} \tau^{-rq^m} (\tau^{rq^m}\sigma^m \tau^s) \tau^{-r} = \tau^{r(1 -  q^m)} \sigma^m \tau^{s},\end{align*}
then, for the cosets of $S$ we have
\begin{align*}\tau^{r} (\sigma^m \tau^s) \tau^{-r} S =   \tau^{r(1 -  q^m)} \sigma^m \tau^{s} S =  \tau^{r(1 -  q^m)} S. \end{align*}
Then $\tau^{r(1 -  q^m)} S = S$ if and only if $r(1-q^m) = \ell t$ for some $\ell$, which implies that ${\displaystyle \frac{t}{r} \, \mid \,1-q^m}$.

Now consider \eqref{st2}. Using \eqref{id-4} twice we obtain
\begin{align*}  \sigma^{\alpha}(\sigma^m \tau^s)\sigma^{-\alpha} =  \sigma^{m}(\sigma^\alpha \tau^s)\sigma^{-\alpha} = \sigma^{m} \tau^{s q^\alpha} = \sigma^m \tau^{sq^\alpha - s} \tau^s = \tau^{s(q^\alpha - 1)q^m }\sigma^m \tau^s,\end{align*}
then for the cosets of $S$ we have
  $$\sigma^{\alpha}(\sigma^m \tau^s)\sigma^{-\alpha}S = \tau^{s(q^\alpha - 1)q^m } S,$$
and $\tau^{s(q^\alpha - 1)q^m } S = S$ if and only if $t \, | \, s(q^\alpha - 1)q^m $. Since by assumption ${\rm gcd}(q,t)=1$, we have that $\tau^{s(q^\alpha - 1)q^m } S = S$ if and only if $t \, | \, s(q^\alpha - 1)$.

Finally, consider \eqref{st3}. Using \eqref{id-4} we obtain
$$ \sigma^{\alpha} \tau^{t} \sigma^{-\alpha} = \tau^{tq^\alpha} \sigma^{\alpha} \sigma^{-\alpha} = \tau^{t q^\alpha},$$
so $\sigma^{\alpha} \tau^{t} \sigma^{-\alpha} S = S$.
\endproof

%%%%%%%%%%%%%%%%%%%%%%%%%%%%%%%%%%%%%%%%%%%%%%%%%%%%%%%

\end{document}